\documentclass[11pt,a4paper]{article}
\usepackage{epsf, graphicx}
\usepackage{latexsym,amsfonts,amsbsy,amssymb}
\usepackage{amsmath,amsthm}
\usepackage{cite}
\usepackage{subfigure}

\makeatletter

\@addtoreset{equation}{section} \makeatother
\newtheorem{theorem}{Theorem}[section]

\makeatletter \@addtoreset{equation}{section} \makeatother
\begin{document}

\title{\bf New metric tensors for anisotropic mesh generation}
\date{}
 \author{Hehu Xie\footnote{LSEC,
ICMSEC, Academy of Mathematics and Systems Science, CAS, Beijing
100080, China \ email: hhxie@lsec.cc.ac.cn} \quad Xiaobo
Yin\footnote{Department of Mathematics, Central China Normal
University, Wuhan 430079, China \ email: yinxb@lsec.cc.ac.cn}}
 \maketitle
 \begin{quote}
\begin{small}
{\bf Abstract.}\,\, A new anisotropic mesh adaptation strategy for
finite element solution of elliptic differential equations is
presented. It generates anisotropic adaptive meshes as quasi-uniform
ones in some metric space, with the metric tensor being computed
based on a posteriori error estimates proposed in \cite{YinXie}. The
new metric tensor explores more comprehensive information of
anisotropy for the true solution than those existing ones. Numerical
results show that this approach can be successfully applied to deal
with poisson and steady convection-dominated problems. The superior
accuracy and efficiency of the new metric tensor to others is
illustrated on various numerical examples of complex two-dimensional
simulations.

{\bf Keywords.} anisotropic; mesh adaptation; metric tensor.

{\bf AMS subject classification.} 65N30, 65N50
\end{small}
\end{quote}
\section{Introduction}
Nowadays, many computational simulations of partial differential
equations (PDEs) involve adaptive triangulations. Mesh adaptation
aims at improving the efficiency and the accuracy of numerical
solutions by concentrating more nodes in regions of large solution
variations than in other regions of the computational domain. As a
consequence, the number of mesh nodes required to achieve a given
accuracy can be dramatically reduced thus resulting in a reduction
of the computational cost. Traditionally, isotropic mesh adaptation
has received much attention, where almost regular mesh elements are
only adjusted in size based on an error estimate. However, in
regions of large solution gradient, adaptive isotropic meshes
usually contain too many elements. Moreover, if the problem at hand
exhibits strong anisotropic features that their solutions change
more significantly in one direction than the others, like boundary
layers, shock waves, interfaces, and edge singularities, etc.. In
such cases it is advantageous to reflect this anisotropy in the
discretization by using meshes with anisotropic elements (sometimes
also called elongated elements). These elements have a small mesh
size in the direction of the rapid variation of the solution and a
larger mesh size in the perpendicular direction. Indeed anisotropic
meshes have been used successfully in many areas, for example in
singular perturbation and flow problems
\cite{Ait-Ali-Yahia,ApelLube,Becker,HabDomBourAitForVal,HabForDomValBou,PerVahMorZ,ZWu}
and in adaptive procedures
\cite{Borouchaki1,Borouchaki2,BuDa,CasHecMohPir,Hecht,PerVahMorZ,Rachowicz}.
For problems with very different length scales in different spatial
directions, long and thin triangles turn out to be better choices
than shape regular ones if they are properly used.

Compared to traditionally used isotropic ones, anisotropic meshes
are more difficult to generate, requiring a full control of both the
shape, size, and orientation of elements. It is necessary to have as
much information as possible on the nature and local behavior of the
solution. We need to convert somehow what we know about the solution
to something having the dimension of a length and containing
directional information. In practice, they are commonly generated as
quasi-uniform meshes in the metric space determined by a tensor (or
a matrix-valued function) specifying the shape, size, and
orientation of elements on the entire physical domain.

Such a metric tensor is often given on a background mesh, either
prescribed by the user or chosen as the mesh from the previous
iteration in an adaptive solver. So far, several meshing strategies
have been developed for generating anisotropic meshes according to a
metric prescription. Examples include blue refinement
\cite{Kornhuber,Lang}, directional refinement \cite{Rachowicz},
Delaunay-type triangulation method
\cite{Borouchaki1,Borouchaki2,CasHecMohPir,PerVahMorZ}, advancing
front method \cite{Garimella}, bubble packing method
\cite{Yamakawa}, local refinement and modification
\cite{HabDomBourAitForVal,Remacle}. On the other hand, variational
methods have received much attention in the recent years typically
for generating structured meshes as they are especially well suited
for finite difference schemes. In these approaches, an adaptive mesh
is considered as the image of a computational mesh under a
coordinate transformation determined by a functional
\cite{Brackbill,Dvinsky,HuangSiam,Jacquotte,Knupp,LiTangZhang}.
Readers are referred to \cite{FreyGeorge} for an overview.

Among these meshing strategies, the definition of the metric tensor
based on the Hessian of the solution seems nowadays commonly
generalized in the meshing community. This choice is largely
motivated by the interesting numerical results obtained by the
results of D'Azevedo \cite{Azevedo1991}, D'Azevedo and Simpson
\cite{AzevedoSimpson} on linear interpolation for quadratic
functions on triangles. For example, Castro-D$\acute{\imath}$az et
al. \cite{CasHecMohPir}, Habashi et al. \cite{HabDomBourAitForVal},
and Remacle et al. \cite{Remacle} define their metric tensor as
\begin{eqnarray}\label{1.1}
\mathcal{M}=|H(u)|\equiv R\left(
\begin{array}{cc}
|\lambda_1| & 0 \\
0 & |\lambda_2|
\end{array}
\right)R^T
\end{eqnarray}
where the Hessian of function $u$ has the eigen-decomposition
$H(u)=R\,\,\mbox{diag}(\lambda_1,\lambda_2)R^T$. To guarantee its
positive definiteness and avoid unrealistic metric, $\mathcal{M}$ in
(\ref{1.1}) is modified by imposing the maximal and minimal edge
lengths. Hecht \cite{Hecht} uses
\begin{eqnarray}\label{1.2_1}
\mathcal{M}=\frac{1}{\epsilon_0\cdot
\mbox{Coef}^2}\frac{|H(u)|}{\max\{\mbox{CutOff},|u|\}}
\end{eqnarray}
for the relative error and
\begin{eqnarray}\label{1.2_2}
\mathcal{M}=\frac{1}{\epsilon_0\cdot
\mbox{Coef}^2}\frac{|H(u)|}{\sup(u)-\inf(u)}
\end{eqnarray}
for the absolute error, where $\epsilon_0$, Coef, and CutOff are the
user specified parameters used for setting the level of the linear
interpolation error (with default value 10$^{-2}$), the value of a
multiplicative coefficient on the mesh size (with default value 1),
and the limit value of the relative error evaluation (with default
value 10$^{-5}$), respectively. In \cite{George}, George and Hecht
define the metric tensor for various norms of the interpolation
error as
\begin{eqnarray}\label{1.3}
\mathcal{M}={\Big (}\frac{c_0}{\epsilon_0}{\Big )}^\nu R\left(
\begin{array}{cc}
|\lambda_1|^\nu & 0 \\
0 & |\lambda_2^\nu
\end{array}
\right)R^T
\end{eqnarray}
where $c_0$ is a constant, $\epsilon_0$ is a given error threshold,
and $\nu=1$ for the $L^{\infty}$ norm and the $H^1$ semi-norm and
$\nu=1/2$ for $L^2$ norm of the error. It is emphasized that the
definitions (\ref{1.1})-(\ref{1.2_2}) are based on the results of
\cite{Azevedo1991} while (\ref{1.3}) largely on heuristic
considerations. Huang \cite{Huang} develop the metric tensors as
\begin{eqnarray}\label{1.4_1}
\mathcal{M}_0=\frac{1}{\sigma}\cdot{\Big
(}\frac{\alpha}{\epsilon_0}{\Big )}\det{\Big
(}I+\frac{1}{\alpha}|H(u)|{\Big )}^{-\frac{1}{6}}{\Big
[}I+\frac{1}{\alpha}|H(u)|{\Big ]}
\end{eqnarray}
for the $L^2$ norm and
\begin{eqnarray}\label{1.4_2}
\mathcal{M}_1=\frac{1}{\sigma}\cdot{\Big
(}\frac{\alpha}{\epsilon_0}{\Big )}^2{\Big
[}I+\frac{1}{\alpha}|H(u)|{\Big ]}
\end{eqnarray}
for the $H^1$ semi-norm.

This list is certainly incomplete, but from the papers we can draw
two conclusions. First, compared with isotropic mesh, significant
improvements in accuracy and efficiency can be gained when a
properly chosen anisotropic mesh is used in the numerical solution
for a large class of problems which exhibit anisotropic solution
features. Second, there are still challenges to access fully
anisotropic information from the computed solution.

The objective of this paper is to develop a new way to get metric
tensors for anisotropic mesh generation in two dimension, which
explores more comprehensive anisotropic information than some exist
methods. The development is based on the error estimates obtained in
our recent work \cite{YinXie} on linear interpolation for quadratic
functions on triangles. These estimates are anisotropic in the sense
that they allow a full control of the shape of elements when used
within a mesh generation strategy.

The application of the error estimates of \cite{YinXie} to formulate
the metric tensor, $\mathcal{M}$, is based on two factors: on the
one hand, as a common practice in the existing anisotropic mesh
generation codes, we assume that the anisotropic mesh is generated
as a quasi-uniform mesh in the metric tensor $\mathcal{M}$, i.e., a
mesh where the elements are equilateral or isosceles right triangle
or other quasi-uniform triangles (shape requirement) in
$\mathcal{M}$ and unitary in size (size requirement) in
$\mathcal{M}$. On the other hand, the anisotropic mesh is required
to minimize the error for a given number of mesh points (or
equidistribute the error). Then $\mathcal{M}$ is constructed from
these requirements. We will establish new metric tensors as
\begin{eqnarray}
\mathcal{M}_0({\bf
x})=\frac{N}{\sigma_0}\det{\Big(}\mathcal{H}{\Big)}^{-\frac{1}{6}}{\Big[}\mathcal{H}{\Big]}
\end{eqnarray}
for the $L^2$ norm and
\begin{eqnarray}
\mathcal{M}_1({\bf x})=\frac{N}{\sigma_1}{\Big
[}\frac{\mbox{tr}(\mathcal{H})}{\sqrt{\det(\mathcal{H})}}{\Big
]}^{\frac{1}{2}}{\Big [}\mathcal{H}{\Big ]}
\end{eqnarray}
for the $H^1$ semi-norm, where $N$ is the number of elements in the
triangulation. Under the condition
$\mathcal{H}=I+\frac{1}{\alpha}|H(u)|$, our metric tensor
$\mathcal{M}_0$ for the $L^2$ norm is similar to (\ref{1.4_1}).
However, the new metric tensor $\mathcal{M}_1$ is essentially
different with those metric tensors mentioned above. The difference
lies on the term
\begin{eqnarray}\label{1.5}
{\Big [}\frac{\mbox{tr}(\mathcal{H})}{\sqrt{\det(\mathcal{H})}}{\Big
]}^{\frac{1}{2}},
\end{eqnarray}
which indicates our metric tensor explores more comprehensive
anisotropic information of the solution when the term (\ref{1.5})
varies significantly at different points or elements. In addition,
numerical results will show that the more anisotropic the solution
is, the more obvious the superiority of the new metric tensor is.

The paper is organized as follows. In Section 2, we describe the
anisotropic error estimates on linear interpolation for quadratic
functions on triangles obtained in the recent work \cite{YinXie}.
The formulation of the metric tensor is developed in Section 3.
Numerical results are presented in Section 4 for various examples of
complex two-dimensional simulations. Finally, conclusions are drawn
in Section 5.

\section{Anisotropic estimates for interpolation error}
Needless to say, the interpolation error depends on the solution and
the size and shape of the elements in the mesh. Understanding this
relation is crucial for the generation of efficient and effective
meshes for the finite element method. In the mesh generation
community, this relation is studied more closely for the model
problem of interpolating quadratic functions. This model is a
reasonable simplification of the cases involving general functions,
since quadratic functions are the leading terms in the local
expansion of the linear interpolation errors. For instance, Nadler
\cite{Nadler} derived an exact expression for the $L^2$-norm of the
linear interpolation error in terms of the three sides ${\bf
\ell}_1$, ${\bf \ell}_2$, and ${\bf \ell}_3$ of the triangle K (see
Figure \ref{Affine_map}),
\begin{eqnarray}\label{Nadler}
\|u-u_I\|^2_{L^2(K)}=\frac{|K|}{180}{\Big[}{\Big(}d_1+d_2+d_3{\Big)}^2+d_1d_2+d_2d_3+d_1d_3{\Big]}.
\end{eqnarray}
where $|K|$ is the area of the triangle, $d_i = {\bf \ell}_i\cdot
H{\bf \ell}_i$ with $H$ being the Hessian of $u$. Bank and Smith
\cite{BankSmith} gave a formula for the $H^1$-seminorm of the linear
interpolation error
\begin{eqnarray*}
\|\nabla(u-u_I)\|^2_{L^2(K)}=\frac{1}{4}{\bf v}\cdot B{\bf v},
\end{eqnarray*}
where ${\bf v}=[d_1,d_2,d_3]^T$,
\begin{eqnarray*}
B=\frac{1}{48|K|}\left(
\begin{array}{ccc}
 |{\bf \ell}_1|^2+|{\bf \ell}_2|^2+|{\bf \ell}_3|^2 & 2{\bf \ell}_1\cdot{\bf \ell}_2 & 2{\bf \ell}_1\cdot{\bf \ell}_3 \\
 2{\bf \ell}_1\cdot{\bf \ell}_2 & |{\bf \ell}_1|^2+|{\bf \ell}_2|^2+|{\bf \ell}_3|^2  & 2{\bf \ell}_2\cdot{\bf \ell}_3\\
 2{\bf \ell}_1\cdot{\bf \ell}_3 & 2{\bf \ell}_2\cdot{\bf \ell}_3 &|{\bf \ell}_1|^2+|{\bf \ell}_2|^2+|{\bf
 \ell}_3|^2
\end{array}
\right).
\end{eqnarray*}
Cao \cite{cao} derived two exact formulas for the $H^1$-seminorm and
$L^2$-norm of the interpolation error in terms of the area, aspect
ratio, orientation, and internal angles of the triangle.

Chen, Sun and Xu \cite{ChenSunXu} obtained the error estimate
\begin{eqnarray*}
\|\nabla(u-u_I)\|^2_{L^p(\Omega)}\leq
CN^{-\frac{2}{n}}\|\sqrt[n]{\det H}\|_{L^{\frac{pn}{2p+n}}(\Omega)},
1\leq p\leq \infty,
\end{eqnarray*}
where the constant $C$ does not depend on $u$ and $N$. They also
showed the estimate is optimal in the sense that it is a lower bound
if $u$ is strictly convex or concave.

Assuming $u=\lambda_1x^2+\lambda_2y^2$, D'Azevedo and Simpson
\cite{Azevedo1989} derived the exact formula for the maximum norm of
the interpolation error
\begin{eqnarray*}
\|(u-u_I)\|^2_{L^{\infty}(K)}=\frac{D_{12}D_{23}D_{31}}{16\lambda_1\lambda_2|K|^2},
\end{eqnarray*}
where $D_{ij}={\bf \ell}_i\cdot \mbox{diag}(\lambda_1,\lambda_2){\bf
\ell}_i$. Based on the geometric interpretation of this formula,
they proved that for a fixed area the optimal triangle, which
produces the smallest maximum interpolation error, is the one
obtained by compressing an equilateral triangle by factors
$\sqrt{\lambda_1}$ and $\sqrt{\lambda_1}$ along the two eigenvectors
of the Hessian of $u$. Furthermore, the optimal incidence for a
given set of interpolation points is the Delaunay triangulation
based on the stretching map (by factors $\sqrt{\lambda_1}$ and
$\sqrt{\lambda_1}$ along the two eigenvector directions) of the grid
points. Rippa \cite{Rippa} showed that the mesh obtained this way is
also optimal for the $L^p$-norm of the error for any $1\leq
p\leq\infty$.

The element-wise error estimates in the following theorem are
developed in \cite{YinXie} using the theory of interpolation for
finite elements and proper numerical quadrature formula.
\begin{theorem}\label{Theorem2.1}
Let $u$ be a quadratic function and $u_I$ is the Lagrangian linear
finite element interpolation of $u$. The following relationship
holds:
\begin{eqnarray*}
\|\nabla(u-u_I)\|^2_{L^2(K)}=\frac{1}{48|K|}\sum_{i=1}^{3}({\bf
\ell}_{i+1}\cdot H{\bf \ell}_{i+2})^2|{\bf \ell}_i|^2,
\end{eqnarray*}
where we prescribe $i+3=i,i-3=i$.
\end{theorem}
Naturally,
\begin{eqnarray*}
\eta_I=\sqrt{\sum_{K\in
\mathcal{T}_h}\frac{1}{48|K|}\sum_{i=1}^{3}({\bf \ell}_{i+1}\cdot
H{\bf \ell}_{i+2})^2|{\bf \ell}_i|^2}
\end{eqnarray*}
is set as the a posteriori estimator for
$\|\nabla(u-u_I)\|_{L^2(\Omega)}$. Numerical experiments in
\cite{YinXie} show that the estimators are always asymptotically
exact under various isotropic and anisotropic meshes.

\section{Metric tensors for anisotropic mesh adaptation}
We now use the results of Theorem \ref{Theorem2.1} to develop a
formula for the metric tensor. As a common practice in anisotropic
mesh generation, we assume that the metric tensor, $\mathcal{M}({\bf
x})$, is used in a meshing strategy in such a way that an
anisotropic mesh is generated as a quasi-uniform mesh in the metric
determined by $\mathcal{M}({\bf x})$. Mathematically, this can be
interpreted as the shape, size and equidistribution requirements as
follows.

{\bf The shape requirement.} The elements of the new mesh,
$\mathcal{T}_h$, are (or are close to being) equilateral in the
metric.

{\bf The size requirement.} The elements of the new mesh
$\mathcal{T}_h$ have a unitary volume in the metric, i.e.,
\begin{eqnarray}\label{3.1}
\int_K\sqrt{\det(\mathcal{M}({\bf x}))}d{\bf x}=1,\quad\forall K\in
\mathcal{T}_h.
\end{eqnarray}

{\bf The equidistribution requirement. } The anisotropic mesh is
required to minimize the error for a given number of mesh points (or
equidistribute the error on every element).

We now state the most significant contributions of this paper.
\subsection{Metric tensor for the $H^1$ semi-norm}
\begin{figure}[ht]
  \centering
  \includegraphics[width=12cm]{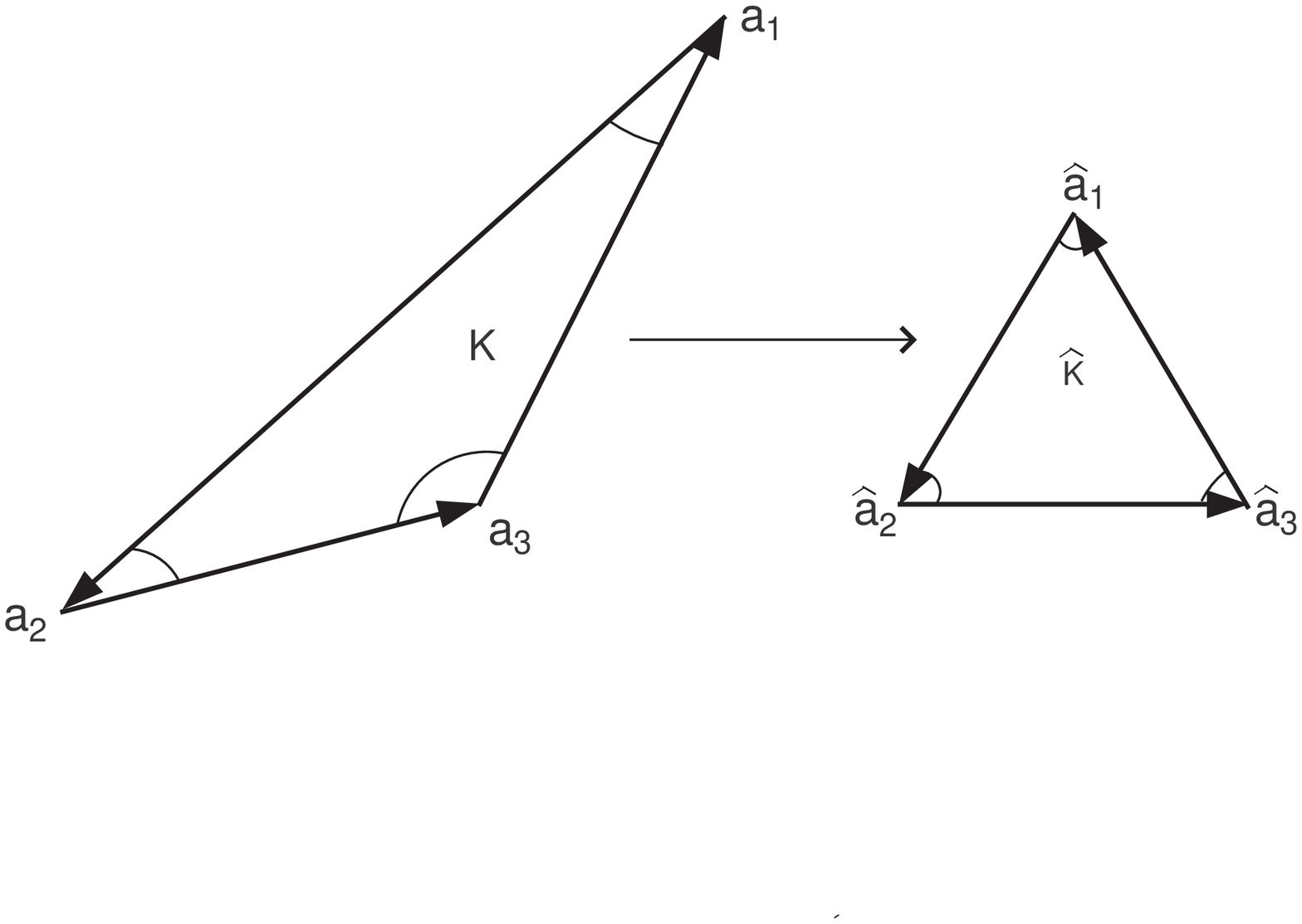}
  \put(-150,181){\large$\mathcal{F}_K$}
  \put(-270,107){\Large$\bf  \ell_1$}
      \put(-178,198){\Large$\bf  \ell_2$}
       \put(-256,183){\Large$\bf  \ell_3$}
  \put(-190,218){\mbox{$\theta_1$}} \put(-292,128){\mbox{$\theta_2$}}
  \put(-240,150){\mbox{$\theta_3$}}
    \put(-73,187){\mbox{$\hat{\theta}_1$}} \put(-100,143){\mbox{$\hat{\theta}_2$}}
  \put(-45,143){\mbox{$\hat{\theta}_3$}}
   \put(-73,123){\mbox{$\hat{\ell}_1$}} \put(-40,165){\mbox{$\hat{\ell}_2$}}
  \put(-98,180){\mbox{$\hat{\ell}_3$}}\vspace{-3.7cm}
  \caption{Affine map ${\bf \hat{x}}=\mathcal{F}_K{\bf x}$ from $K$ to the
reference triangle $\hat{K}$.}\label{Affine_map}
  \end{figure}

Assume that $H({\bf x})$ is a symmetric positive definite matrix on
every point ${\bf x}$. Let $\mathcal{M}_1({\bf x})=\theta_1M_1({\bf
x})$ where $\theta_1$ is to be determined. Here $M_1({\bf x})$ is
often called the monitor function. Both the monitor function
$M_1({\bf x})$ and metric tensor $\mathcal{M}_1({\bf x})$ play the
same role in mesh generation, i.e., they are used to specify the
size, shape, and orientation of mesh elements throughout the
physical domain. The only difference lies in the way they specify
the size of elements. Indeed, $M_1({\bf x})$ specifies the element
size through the equidistribution condition, while
$\mathcal{M}_1({\bf x})$ determines the element size through the
unitary volume requirement (\ref{3.1}).

Assume that $H({\bf x})$ is a constant matrix on $K$, denoted by
$H_K$, then so does $M_1({\bf x})$, denoted by $M_{1,K}$. Since
$H_K$ is a symmetric positive definite matrix, we consider the
singular value decomposition $H_K=R^T\Lambda R$, where
$\Lambda=\mbox{diag}(\lambda_1,\lambda_2)$ is the diagonal matrix of
the corresponding eigenvalues ($\lambda_1,\lambda_2>0$) and $R$ is
the orthogonal matrix having as rows the eigenvectors of $H_K$.
Denote by $F_K$ and ${\bf t}_K$ the matrix and the vector defining
the invertible affine map $\hat{\bf x}=\mathcal{F}_K({\bf
x})=F_K{\bf x} + {\bf t}_K$ from the generic element $K$ to the
reference triangle $\hat{K}$ (see Figure \ref{Affine_map}).

Let $M_{1,K}=C_KH_K$, then
$F_K=C_K^{\frac{1}{2}}\Lambda^{\frac{1}{2}} R$ and
$M_{1,K}=F_K^TF_K$. Mathematically, the shape requirement can be
expressed as
\begin{eqnarray*}
|\hat{\ell}_i|=L\,\,\mbox{and}\,\,
\cos\hat{\theta}_i=\frac{\hat{\ell}_{i+1}\cdot\hat{\ell}_{i+2}}{L^2}=\frac{1}{2},\,i=1,2,3,
\end{eqnarray*}
where $L$ is a constant for every element $K$. Together with Theorem
\ref{Theorem2.1} we have
\begin{eqnarray*}
\|\nabla(u-u_I)\|^2_{L^2(K)}&=&\frac{1}{48|K|}\sum_{i=1}^{3}({\bf
\ell}_{i+1}\cdot H_K{\bf \ell}_{i+2})^2|{\bf
\ell}_i|^2\nonumber\\
&=&\frac{1}{48|K|C_K^2}\sum_{i=1}^{3}({\bf \ell}_{i+1}\cdot
M_{1,K}{\bf
\ell}_{i+2})^2|{\bf \ell}_i|^2\nonumber\\
&=&\frac{L^4}{48|K|C_K^2}\sum_{i=1}^{3}(\cos\hat{\theta}_i)^2|{\bf
\ell}_i|^2\nonumber\\
&=&\frac{C_K\sqrt{\lambda_1\lambda_2}L^4}{192|\hat{K}|C_K^2}\sum_{i=1}^{3}{\Big|}C_K^{-\frac{1}{2}}R^{-1}\Lambda^{-\frac{1}{2}}\hat{\ell}_i{\Big|}^2\nonumber\\
&=&\frac{\sqrt{\lambda_1\lambda_2}L^2}{48\sqrt{3}C_K^2}\sum_{i=1}^{3}{\Big|}\Lambda^{-\frac{1}{2}}\hat{\ell}_i{\Big|}^2\nonumber\\
&=&\frac{L^4}{32\sqrt{3}C_K^2}{\Big(}\sqrt{{\lambda_1}/{\lambda_2}}+\sqrt{{\lambda_2}/{\lambda_1}}{\Big)}.
\end{eqnarray*}
To satisfy the equidistribution requirement, let
\begin{eqnarray*}
\|\nabla(u-u_I)\|^2_{L^2(K)}={\Big (}\sum\limits_{K\in
\mathcal{T}_h}e_K^2{\Big )}/N=\epsilon_1^2/N,
\end{eqnarray*}
where $N$ is the number of elements of $\mathcal{T}_h$. Then
\begin{eqnarray*}
C_K=L^2{\Big (}\frac{N}{32\sqrt{3}\epsilon_1^2}{\Big
)}^{\frac{1}{2}}{\Big
(}\sqrt{{\lambda_1}/{\lambda_2}}+\sqrt{{\lambda_2}/{\lambda_1}}{\Big
)}^{\frac{1}{2}}=\overline{C}{\Big
(}\sqrt{{\lambda_1}/{\lambda_2}}+\sqrt{{\lambda_2}/{\lambda_1}}{\Big
)}^{\frac{1}{2}}=\overline{C}{\Big
[}\frac{\mbox{tr}(H)}{\sqrt{\mbox{det}(H)}}{\Big ]}^{\frac{1}{2}},
\end{eqnarray*}
where $\overline{C}$ is a constant which depends on the prescribed
error and the numbers of elements. Generally speaking, $H({\bf x})$
can not be a constant matrix on $K$ , $M_1({\bf x})$ should be the
form
\begin{eqnarray*}
M_1({\bf x})={\Big [}\frac{\mbox{tr}(H)}{\sqrt{\mbox{det}(H)}}{\Big
]}^{\frac{1}{2}}{\Big [}H(u){\Big ]}={\Big (}\sqrt{{\lambda_1({\bf
x})}/{\lambda_2({\bf x})}}+\sqrt{{\lambda_2({\bf
x})}/{\lambda_1({\bf x})}}{\Big )}^{\frac{1}{2}}{\Big [}H(u){\Big
]},
\end{eqnarray*}
since $M_1({\bf x})$ can be modified by multiplying a constant. Note
that $\lambda_1({\bf x})$ and $\lambda_2({\bf x})$ are the
corresponding eigenvalue of $H(u)$ at point ${\bf x}$.

To establish $\mathcal{M}_1({\bf x})$, the size requirement
(\ref{3.1}) should be used, which leads to
\begin{eqnarray*}
\theta_1\int_K\rho_1({\bf x}) d{\bf x}=1,
\end{eqnarray*}
where
\begin{eqnarray*}
\rho_1({\bf x})=\sqrt{\det(M_1({\bf x}))}.
\end{eqnarray*}
Summing the above equation over all the elements of $\mathcal{T}_h$,
one gets
\begin{eqnarray*}
\theta_1\sigma_1=N,
\end{eqnarray*}
where
\begin{eqnarray*}
\sigma_1=\int_{\Omega}\rho_1({\bf x}) d{\bf x}.
\end{eqnarray*}
Thus, we get
\begin{eqnarray*}
\theta_1=\frac{N}{\sigma_1},
\end{eqnarray*}
and as a consequence,
\begin{eqnarray*}
\mathcal{M}_1({\bf x})=\frac{N}{\sigma_1}{\Big
[}\frac{\mbox{tr}(H)}{\sqrt{\mbox{det}(H)}}{\Big
]}^{\frac{1}{2}}{\Big [}H(u){\Big ]}=\frac{N}{\sigma_1}{\Big
(}\sqrt{{\lambda_1({\bf x})}/{\lambda_2({\bf
x})}}+\sqrt{{\lambda_2({\bf x})}/{\lambda_1({\bf x})}}{\Big
)}^{\frac{1}{2}}{\Big [}H(u){\Big ]}.
\end{eqnarray*}

\subsection{Metric tensor for the $L^2$ norm}
Using the expression (\ref{Nadler}) for the $L^2$-norm of the linear
interpolation error derived by Nadler\cite{Nadler} and the condition
(\ref{3.1}), we have
\begin{eqnarray*}
\|u-u_I\|^2_{L^2(K)}&=&\frac{|K|}{180}{\Big [}{\Big(}\sum_{i=1}^{3}d_i{\Big)}^2+d_1d_2+d_2d_3+d_1d_3{\Big]}\nonumber\\
&=&\frac{|K|}{180C_K^2}{\Big[}{\Big(}\sum_{i=1}^{3}|\hat{\ell}_i|^2{\Big)}^2+\sum_{i=1}^{3}{\Big(}|\hat{\ell}_{i+1}||\hat{\ell}_{i+2}|{\Big)}^2{\Big]}\nonumber\\
&=&\frac{L^4|K|}{15C_K^2}=\frac{L^4|\hat{K}|}{15C_K^3\sqrt{\lambda_1\lambda_2}}=\frac{\sqrt{3}L^6}{60C_K^3\sqrt{\lambda_1\lambda_2}}.
\end{eqnarray*}
To satisfy the equidistribution requirement, let
\begin{eqnarray*}
\|u-u_I\|^2_{L^2(K)}={\Big (}\sum\limits_{K\in
\mathcal{T}_h}e_K^2{\Big )}/N=\epsilon_0^2/N.
\end{eqnarray*}
Using similar argument with last subsection, we easily get the
monitor function
\begin{eqnarray*}
M_0({\bf x})=\det{\Big (}H{\Big )}^{-\frac{1}{6}}{\Big [}H(u){\Big
]},
\end{eqnarray*}
and the metric tensor
\begin{eqnarray*}
\mathcal{M}_0({\bf x})=\frac{N}{\sigma_0}\det{\Big (}H{\Big
)}^{-\frac{1}{6}}{\Big [}H(u){\Big ]}
\end{eqnarray*}
for the $L^2$ norm.
\subsection{Practice use of metric tensors}
So far we assume that $H({\bf x})$ is a symmetric positive definite
matrix at every point. However this assumption doesn't hold in many
cases. In order to obtain a symmetric positive definite matrix, the
following procedure are often implemented. First, the Hessian $H$ is
modified into $|H|=R^T\,\,\mbox{diag}(|\lambda_1|,|\lambda_2|)R$ by
taking the absolute value of its eigenvalues
(\cite{HabForDomValBou}). Since $|H|$ is only semi-positive
definite, $\mathcal{M}_0$ and $\mathcal{M}_1$ cannot be directly
applied to generate the anisotropic meshes. To avoid this
difficulty, we regularize the expression with two flooring
parameters $\alpha_0>0$ for $\mathcal{M}_0$ and $\alpha_1>0$ for
$\mathcal{M}_1$, respectively (\cite{HuangSiam}). Replace $|H|$ with
\begin{eqnarray*}
\mathcal{H}=\alpha_1I+|H|,
\end{eqnarray*}
then we get the modified metric tensor
\begin{eqnarray}\label{H1metric}
\mathcal{M}_1({\bf x})=\frac{N}{\sigma_1}{\Big
[}\frac{\mbox{tr}(\mathcal{H})}{\sqrt{\det(\mathcal{H})}}{\Big
]}^{\frac{1}{2}}{\Big [}\mathcal{H}{\Big ]}.
\end{eqnarray}
Similarly, replacing $|H|$ with
\begin{eqnarray*}
\mathcal{H}=\alpha_0I+|H|
\end{eqnarray*}
leads to
\begin{eqnarray}\label{L2metric}
\mathcal{M}_0({\bf x})=\frac{N}{\sigma_0}\det{\Big
(}\mathcal{H}{\Big )}^{-\frac{1}{6}}{\Big [}\mathcal{H}{\Big ]}.
\end{eqnarray}
The two modified metric tensors (\ref{H1metric}) and
(\ref{L2metric}) are suitable for practical mesh generation.

\section{Numerical experiments}
We have elaborated that our metric tensor explores more anisotropic
information for a given problem. Then it is crucial to check that
whether our metric tensor exhibits better than those without the
term
\begin{eqnarray}\label{factor}
{\Big
[}\frac{\mbox{tr}(\mathcal{H})}{\sqrt{\mbox{det}(\mathcal{H})}}{\Big
]}^{\frac{1}{2}} \quad \mbox{or}\quad {\Big (}\sqrt{{\lambda_1({\bf
x})}/{\lambda_2({\bf x})}}+\sqrt{{\lambda_2({\bf
x})}/{\lambda_1({\bf x})}}{\Big )}^{\frac{1}{2}}.
\end{eqnarray}

\subsection{Mesh adaptation tool}
 All the presented experiments are performed by using the BAMG
software \cite{Hecht}. Given a ¡°background¡± mesh and an
approximation solution, BAMG generates the mesh according to the
metric. The code allows the user to supply his/her own metric tensor
defined on a background mesh. In our computation, the background
mesh has been taken as the most recent mesh available.

\subsection{Comparisons between two types of metric tensors}
\begin{figure}[ht!]
  \vspace{-2cm}\hspace{1cm}
\includegraphics[width=16cm]{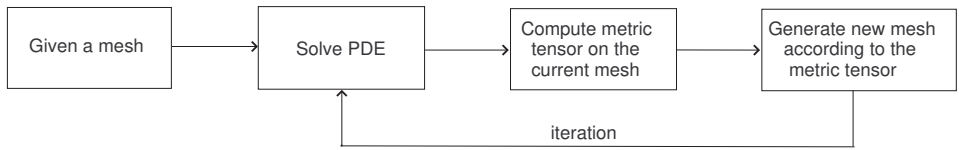}
  \caption{Iterative procedure for adaptive mesh solution of
  PDEs.}\label{adaptive_loop}
\end{figure}

Specifically, in every example, the PDE is discretized using linear
triangle finite elements. Two serials of adaptive meshes of almost
the same number of elements are shown (the iterative procedure for
solving PDEs is shown in Figure \ref{adaptive_loop}).
The finite element solutions are computed by using the metric tensor
of modified Hessian
\begin{eqnarray}\label{modifiedhessian}
\mathcal{M}_1^m({\bf x})=\frac{N}{\sigma_1}|H|,
\end{eqnarray}
and the new metric tensor
\begin{eqnarray}\label{newmetric}
\mathcal{M}_1^n({\bf x})=\frac{N}{\sigma_1}{\Big
[}\frac{\mbox{tr}(|H|)}{\sqrt{\det(|H|)}}{\Big ]}^{\frac{1}{2}}|H|.
\end{eqnarray}
\begin{figure}[ht]
 \includegraphics[width=8cm]{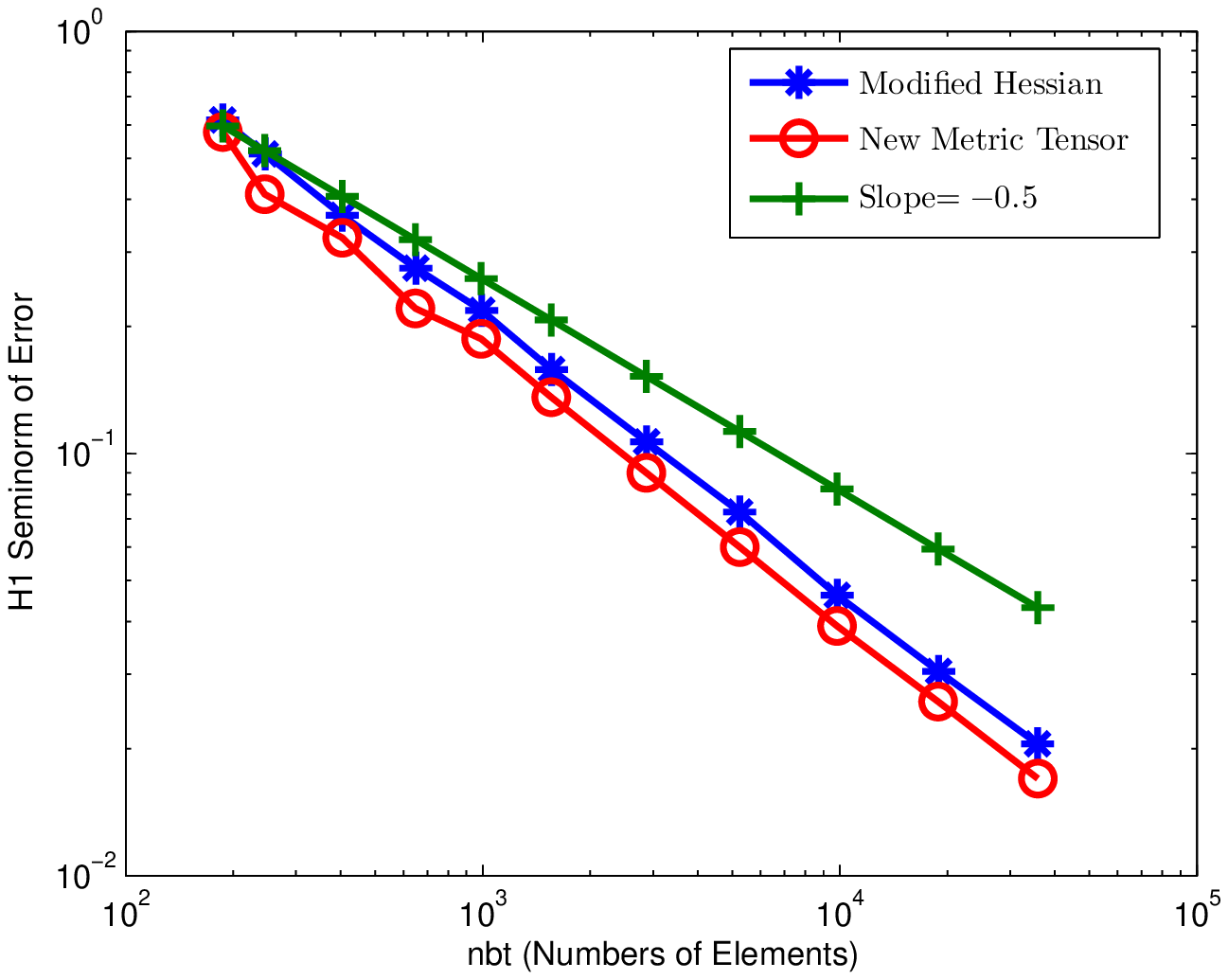}
  \put(-120,167){(a)}
    \put(-15,0){\resizebox{8cm}{!}{\epsffile{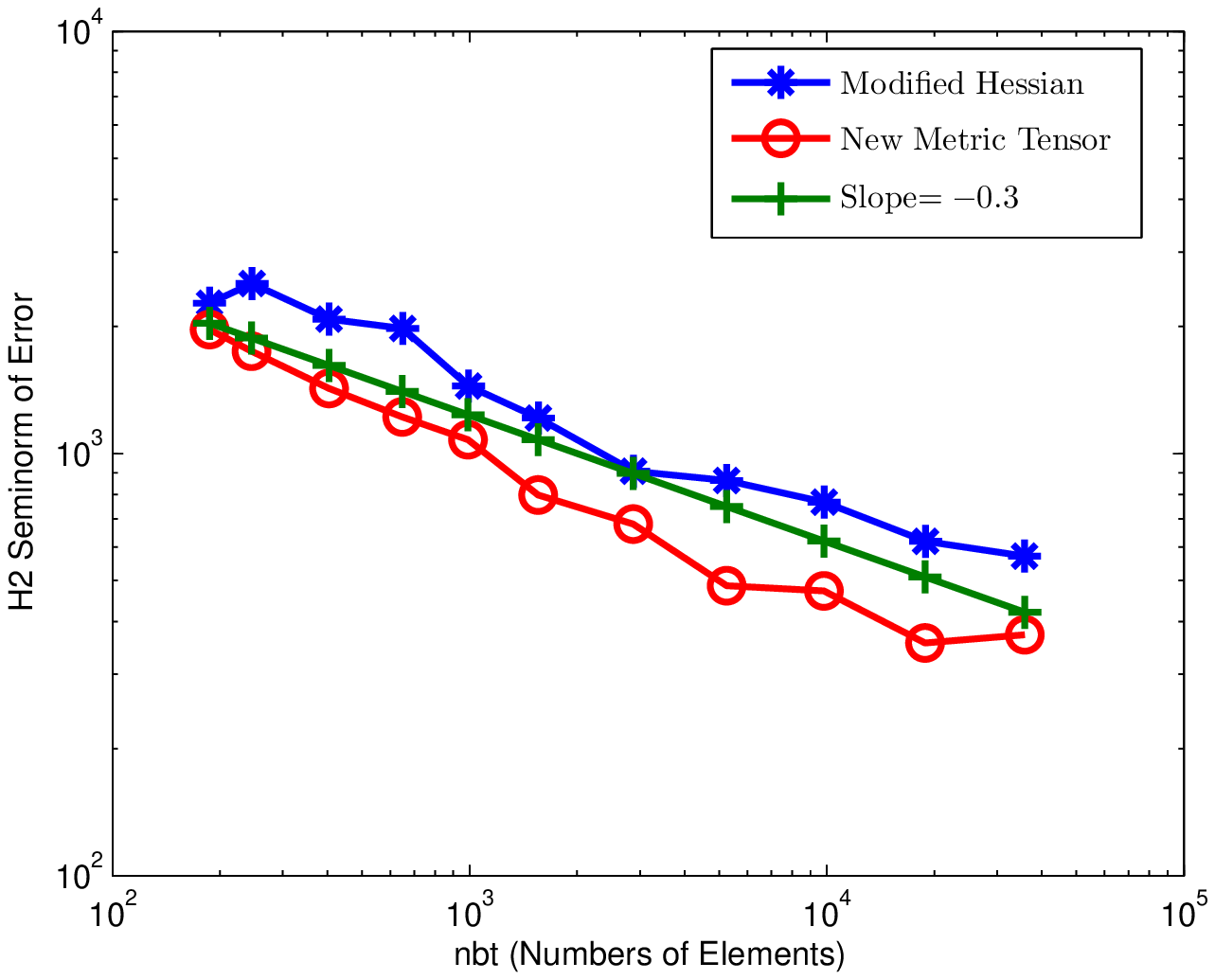}}}
    \put(100,167){(b)}
    \caption{Example 4.1. The $H^1$ and $H^2$ semi-norms of discretization error are plotted
  as functions of the number of elements ($nbt$) in (a) and (b),
  respectively.}\label{Example1_error}
\end{figure}
Notice that the formulas of the metric tensors involve second order
derivatives of the physical solution. Generally speaking, one can
assume that the nodal values of the solution or their approximations
are available, as typically in the numerical solution of partial
differential equations. Then, one can get approximations for the
second order derivatives using gradient recovery techniques (such as
\cite{ZZ1987} and \cite{Zhangzm}) twice or a Hessian recovery
technique (using piecewise quadratic polynomials fitting in
least-squares sense to nodal values of the computed solution
\cite{Zhangxd}) just once, although their convergence has been
analyzed only on isotropic meshes. In our computations, we use the
technique \cite{ZZ1987} twice.

In the current computation, each run is stopped after 10 iterations
to guarantee that the adaptive procedure tends towards stability.

Denote by $nbt$ the number of the elements (triangles in 2D case) in
the current mesh. The number of triangles or nodes is adjusted when
necessary by trial and errors through the modification of the
multiplicative coefficient of the metric tensors.

\noindent{\bf Example 4.1.}  This example, though numerically solved
in $ \Omega\equiv(0,1)\times(0,1),$
 is in fact one-dimensional:
\begin{eqnarray*}
-\kappa\triangle u+\frac{\partial u}{\partial x_1}&=&f,
\end{eqnarray*}
with $f=0$, $u(x_1=0,x_2)=0$, $u(x_1=1,x_2)=1$, and $\frac{\partial
u}{\partial n}=0$ along the top and bottom sides of the square
(taken from \cite{BuDa}). The exact solution is given by
\begin{eqnarray*}
u({\bf x})=\frac{1-e^{\frac{x_1}{\kappa}}}{1-e^{\frac{1}{\kappa}}},
\end{eqnarray*}
with a boundary layer of width $\kappa$ at $x_1\approx1$. We have
set $\kappa=0.0015$, so that convection is dominant and the Galerkin
method yields oscillatory solutions unless the mesh is highly
refined at the boundary layer.

Figure \ref{Example1_error} (a) compares $H^1$ semi-norms of the
discretization errors for the two metric tensors. In (b) $H^2$
semi-norm of error is computed by the difference between the Hessian
of $u$ and the recovered one, which exhibits the quality of meshes
to some certain extent.
\begin{figure}[ht]
\includegraphics[width=8cm]{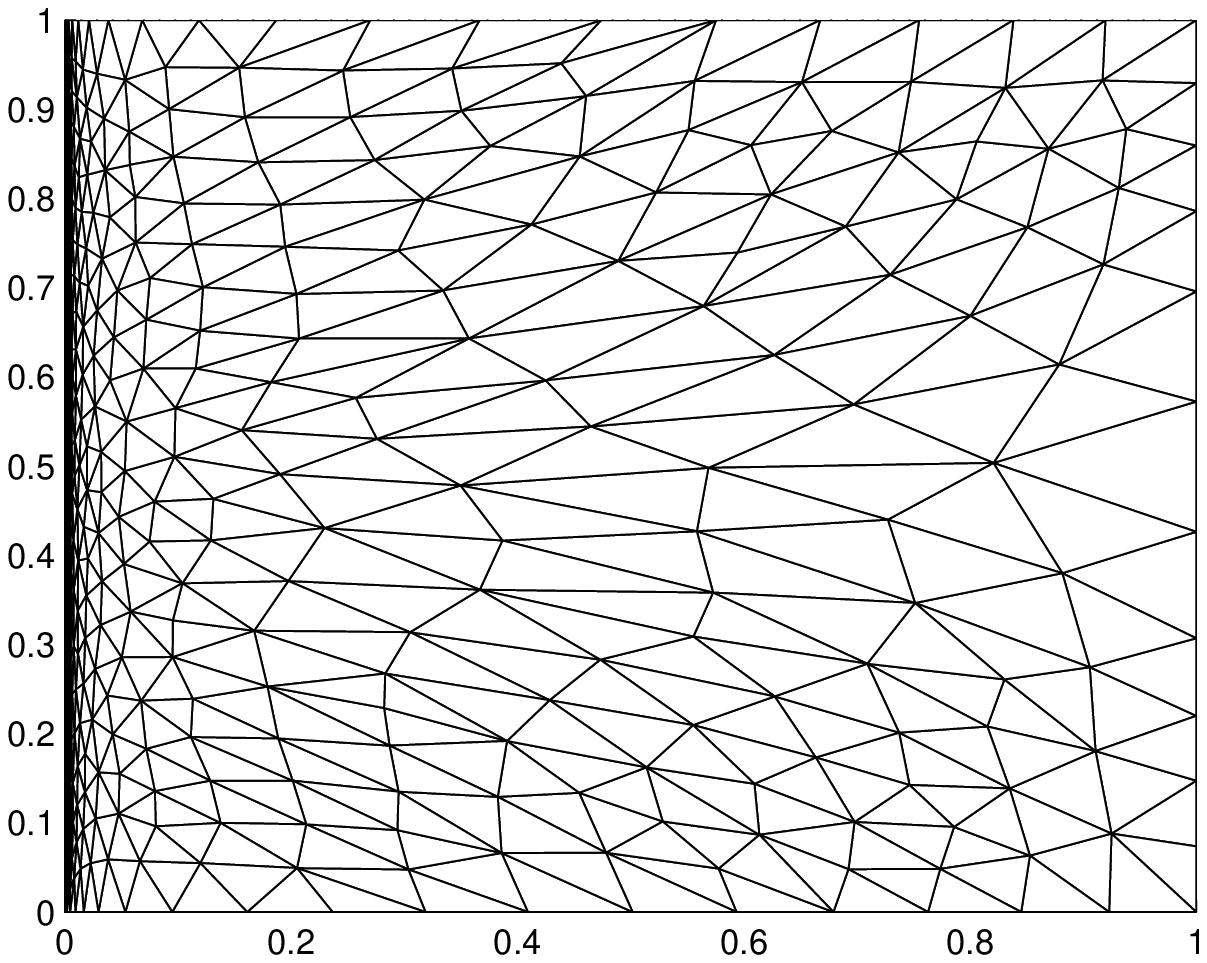}
  \put(-115,167){(a)}
  \put(-20,0){\resizebox{8cm}{!}{\epsffile{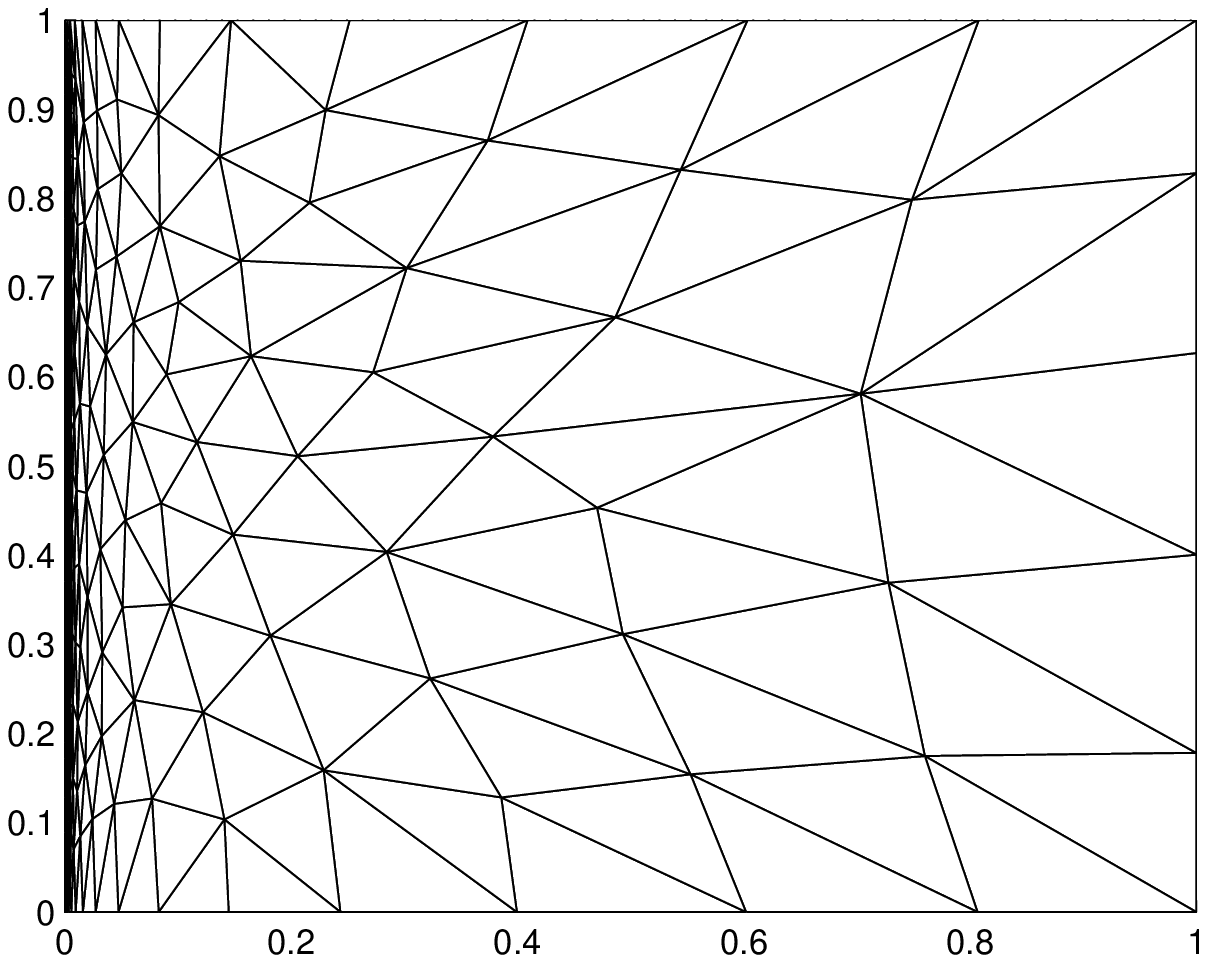}}}
  \put(97,167){(b)}
  \caption{Example 4.2. (a) Anisotropic mesh obtained with the metric tensor $\mathcal{M}_1^m({\bf x})$:
$nbt=4244$, $|e|_{H^1}=0.3727$, and $|e|_{H^2}=1762$. (b)
Anisotropic mesh obtained the metric tensor $\mathcal{M}_1^n({\bf
x})$: $nbt=4243$, $|e|_{H^1}=0.2842$, and
$|e|_{H^2}=1101$.}\label{Example2_mesh}
\end{figure}
\noindent{\bf Example 4.2.} We study the Poisson equation (taken
from \cite{ForPer})
\begin{equation*}
     \left \{
     \begin{array}{lll}
     -\triangle u&=&f,\quad \mbox{in}\quad\Omega\equiv(0,1)\times(0,1),\\
     u&=&0,\quad \mbox{on} \quad\partial \Omega,
     \end{array}
     \right .
\end{equation*}
where $f$ has been chosen in such a way that the exact solution
$u({\bf x})=[1-e^{-\alpha x_1}-(1-e^{-\alpha})x_1]4x_2(1-x_2)$ and
$\alpha$ is chosen to be 1000. Notice that the function $u$ exhibits
an exponential layer along the boundary $x_1=0$ with an initial step
of $\alpha$.
\begin{figure}[ht]
\includegraphics[width=8cm]{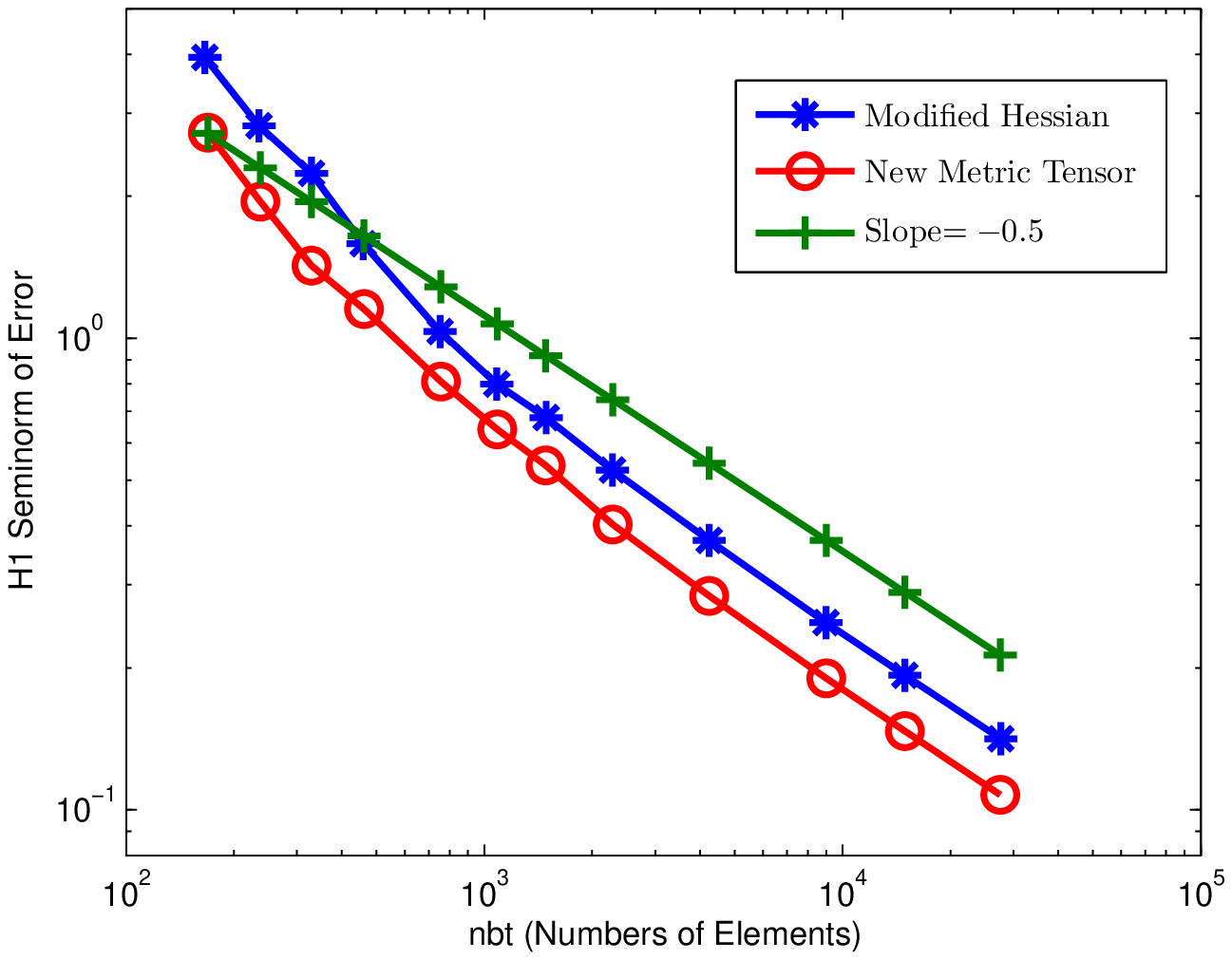}
  \put(-120,167){(a)}
    \put(-15,0){\resizebox{8cm}{!}{\epsffile{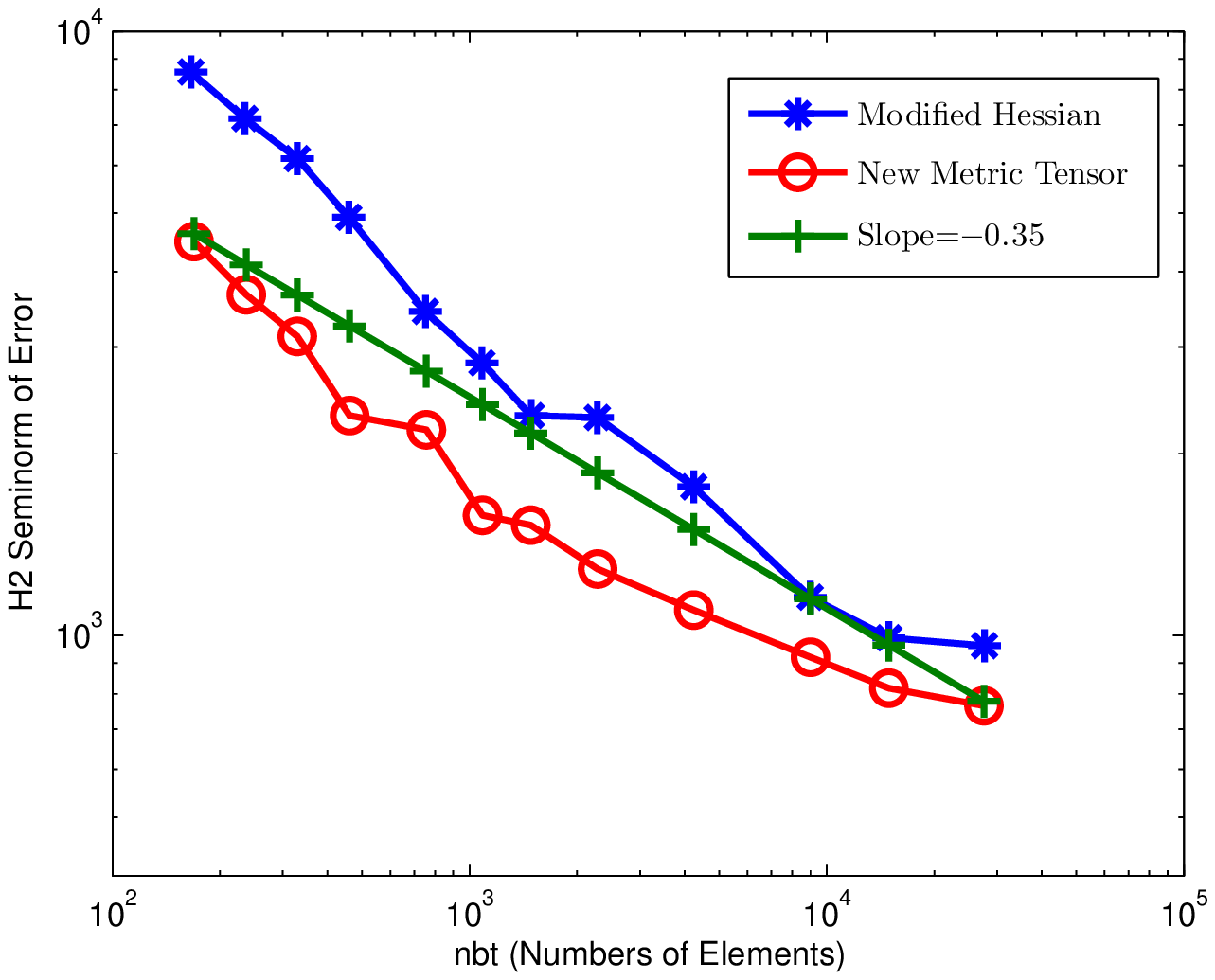}}}
    \put(100,167){(b)}
       \caption{Example 4.2. The $H^1$ and $H^2$ semi-norms of discretization error are plotted
as functions of the number of elements ($nbt$) in (a) and (b),
respectively.}\label{Example2_error}
\end{figure}
\begin{figure}[ht]
\includegraphics[width=8cm]{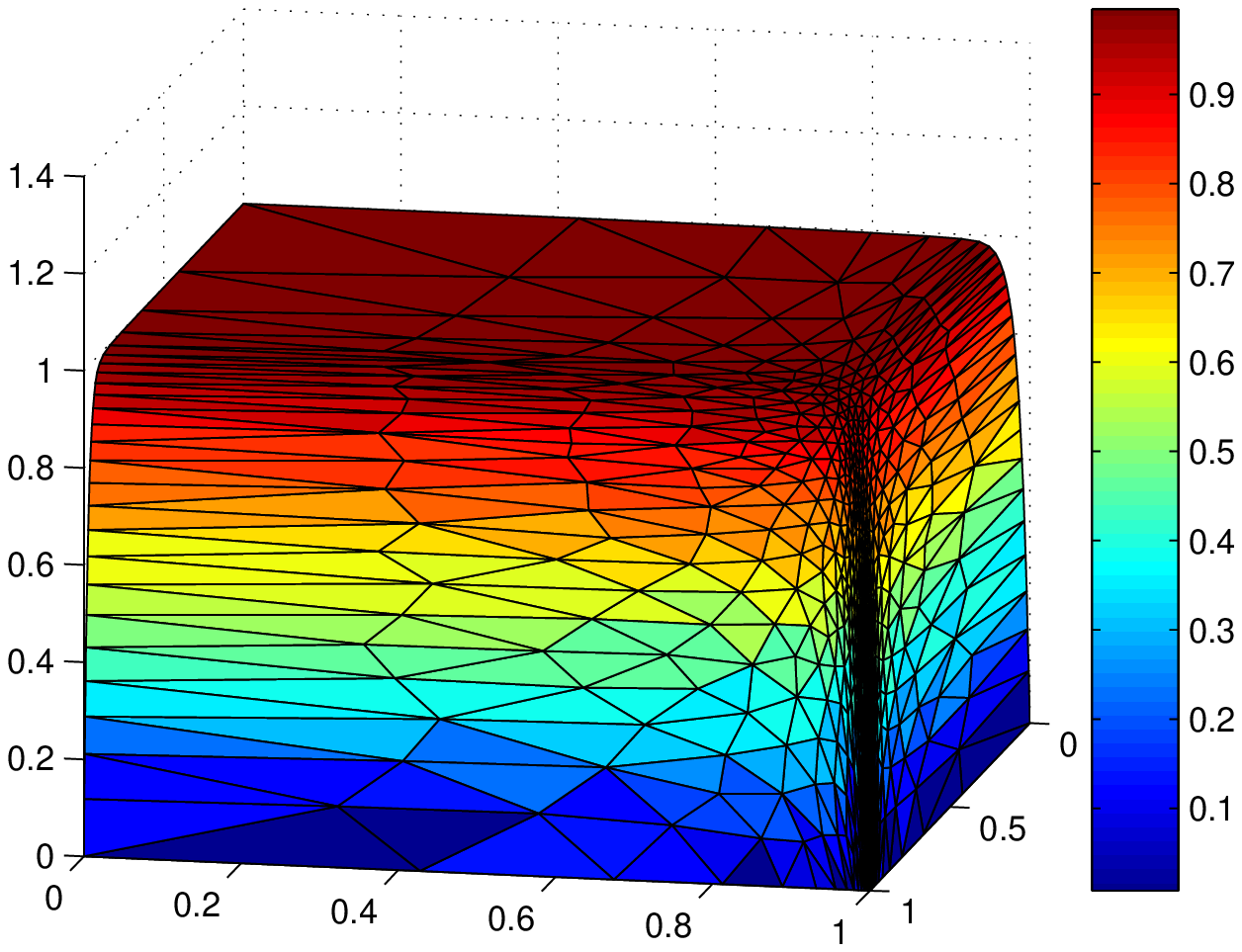}
  \put(-125,170){(a)}
  \put(-10,0){\resizebox{8cm}{!}{\epsffile{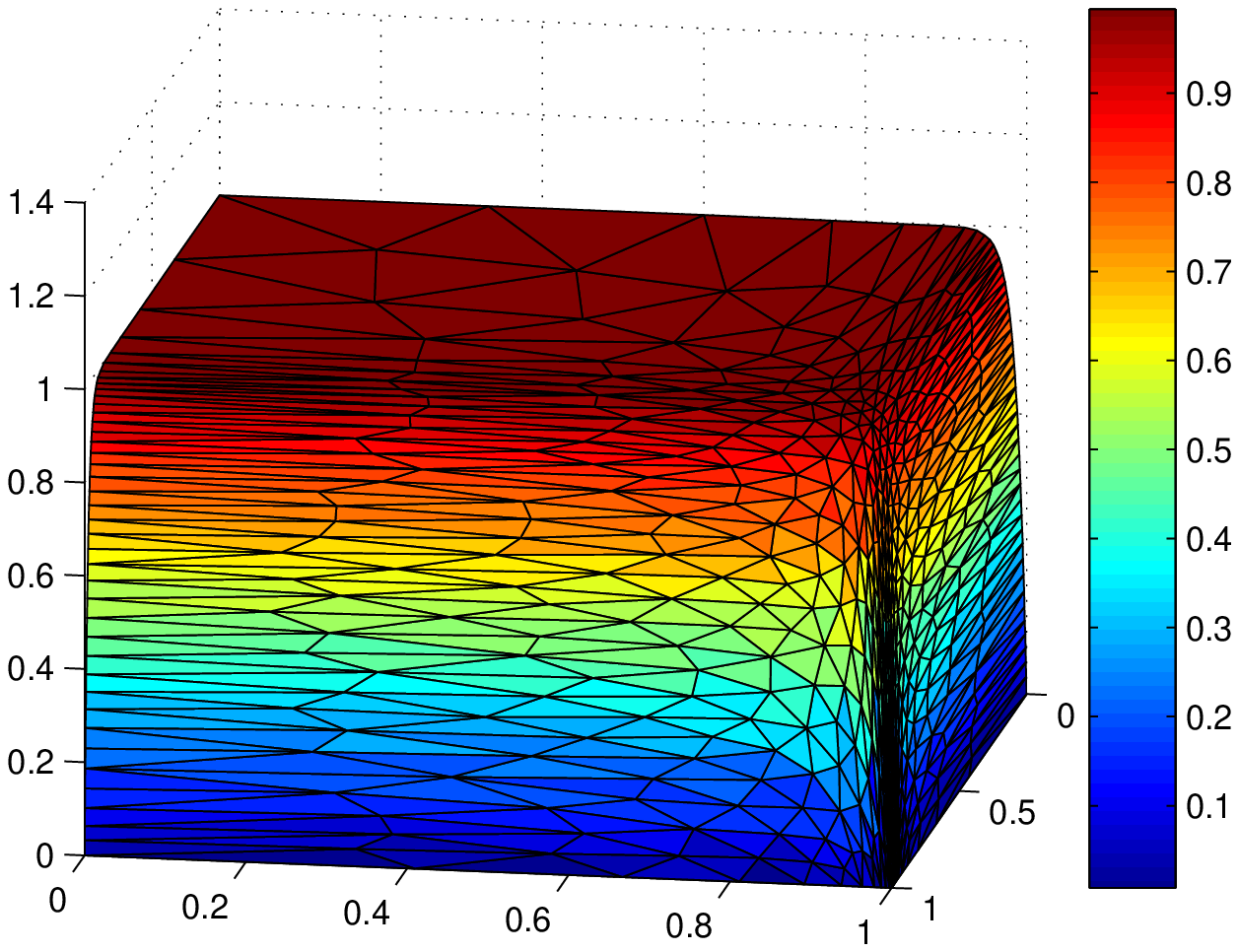}}}
  \put(100,170){(b)}
\caption{Example 4.3. $\beta=40$, (a) Anisotropic mesh obtained with
the metric tensor $\mathcal{M}_1^m({\bf x})$: $nbt=892$,
$|e|_{H^1}=0.2581$, $|e|_{H^2}=102.0$. (b) Anisotropic mesh obtained
with the metric tensor $\mathcal{M}_1^n({\bf x})$: $nbt=891$,
$|e|_{H^1}=0.1893$, $|e|_{H^2}=57.57$.}\label{Example3_solution}
\end{figure}
Figure \ref{Example2_mesh} contains two meshes obtained by the two
different metric tensors (\ref{modifiedhessian}) and
(\ref{newmetric}). It is easily seen that the two meshes are obvious
different. The former explores more anisotropic features of the
solution $u$ than the latter. In other words, the term
(\ref{factor}) complements more comprehensive
 information of the exact solution. Figure \ref{Example2_error} contains $H^1$
and $H^2$ semi-norms of error similar to Figure
\ref{Example1_error}.

\noindent{\bf Example 4.3.} Let $\Omega\equiv(0,1)\times(0,1)$, and
$u$ be the solution of
\begin{eqnarray*}
-\triangle
u=\beta(\beta-1)x_1^{\beta-2}(1-x_2^{2\beta})+2\beta(2\beta-1)x_2^{2\beta-2}(1-x_1^{\beta})
\end{eqnarray*}
with boundary conditions $u=0$ along the sides $x_1=1$ and $x_2=1$,
and $\partial u/\partial n=0$ along $x_1=0$ and $x_2=0$ (taken from
\cite{BuDa}). The exact solution $u({\bf
x})=(1-x_1^{\beta})(1-x_2^{2\beta})$ exhibits two boundary layers
along the right and top sides, the latter being stronger than the
former.

\begin{figure}[ht]
\includegraphics[width=8cm]{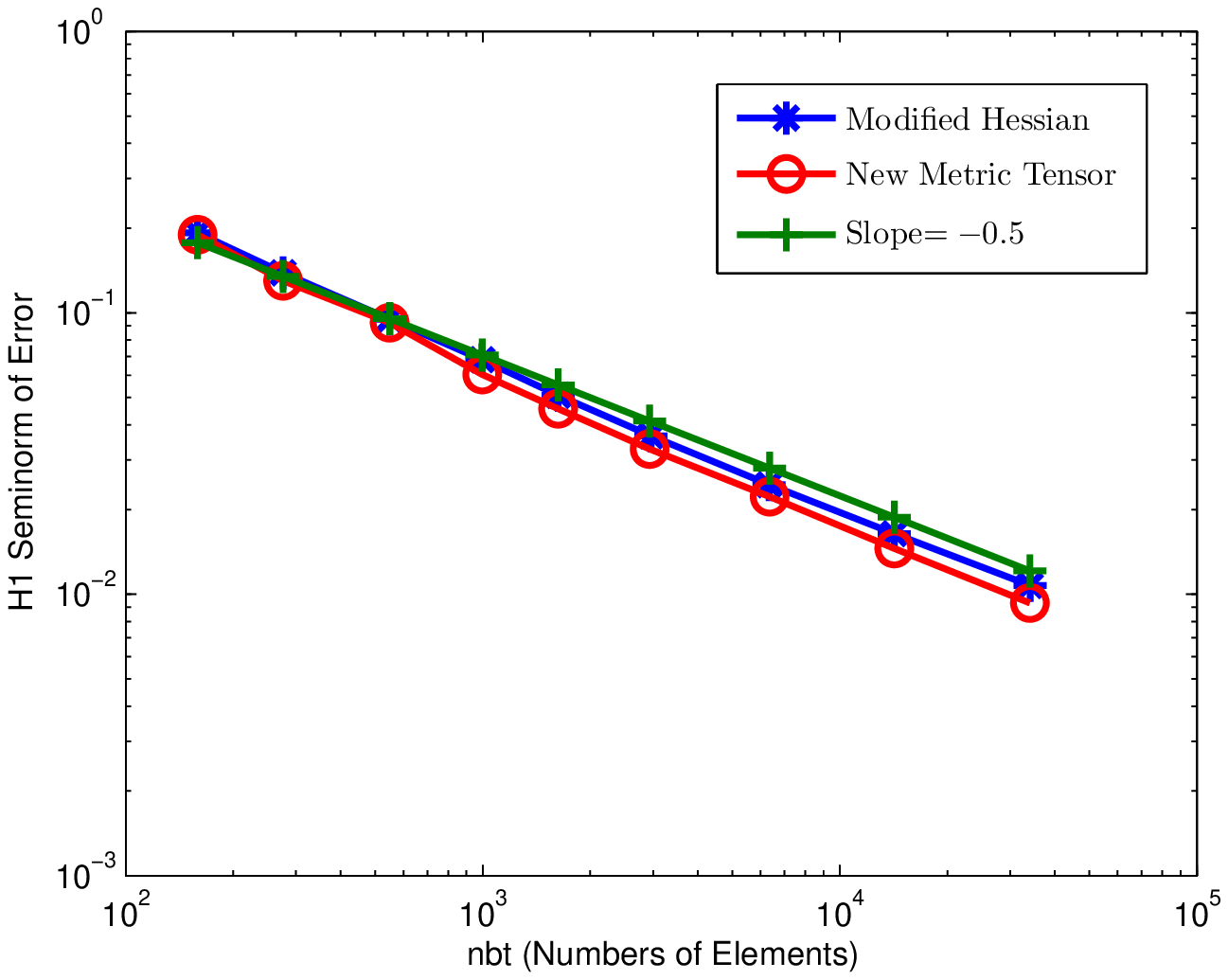}
   \put(-125,167){(a)}
  \put(-20,0){\resizebox{8cm}{!}{\epsffile{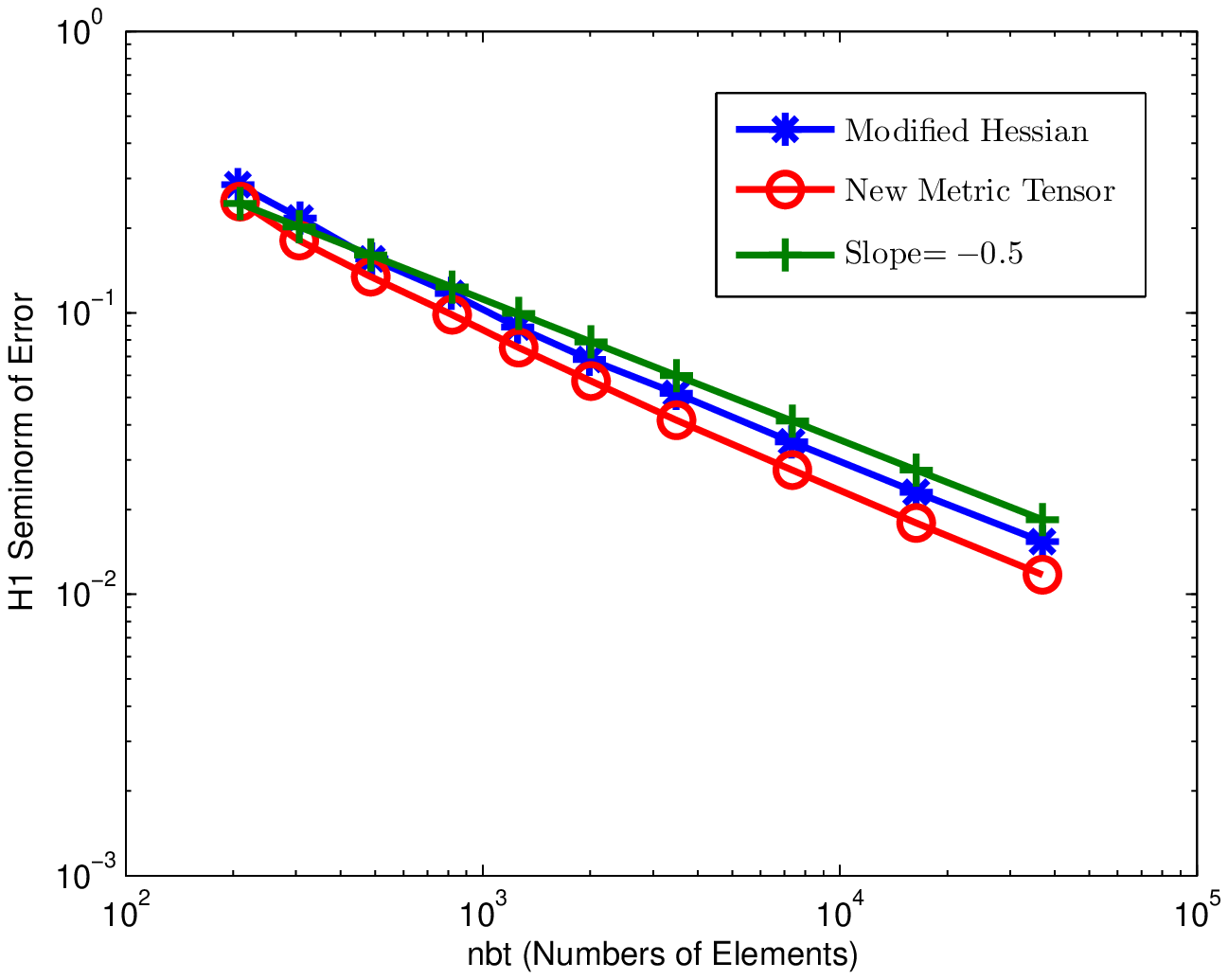}}}
  \put(100,167){(b)}
  \put(-228,-180){\resizebox{8cm}{!}{\epsffile{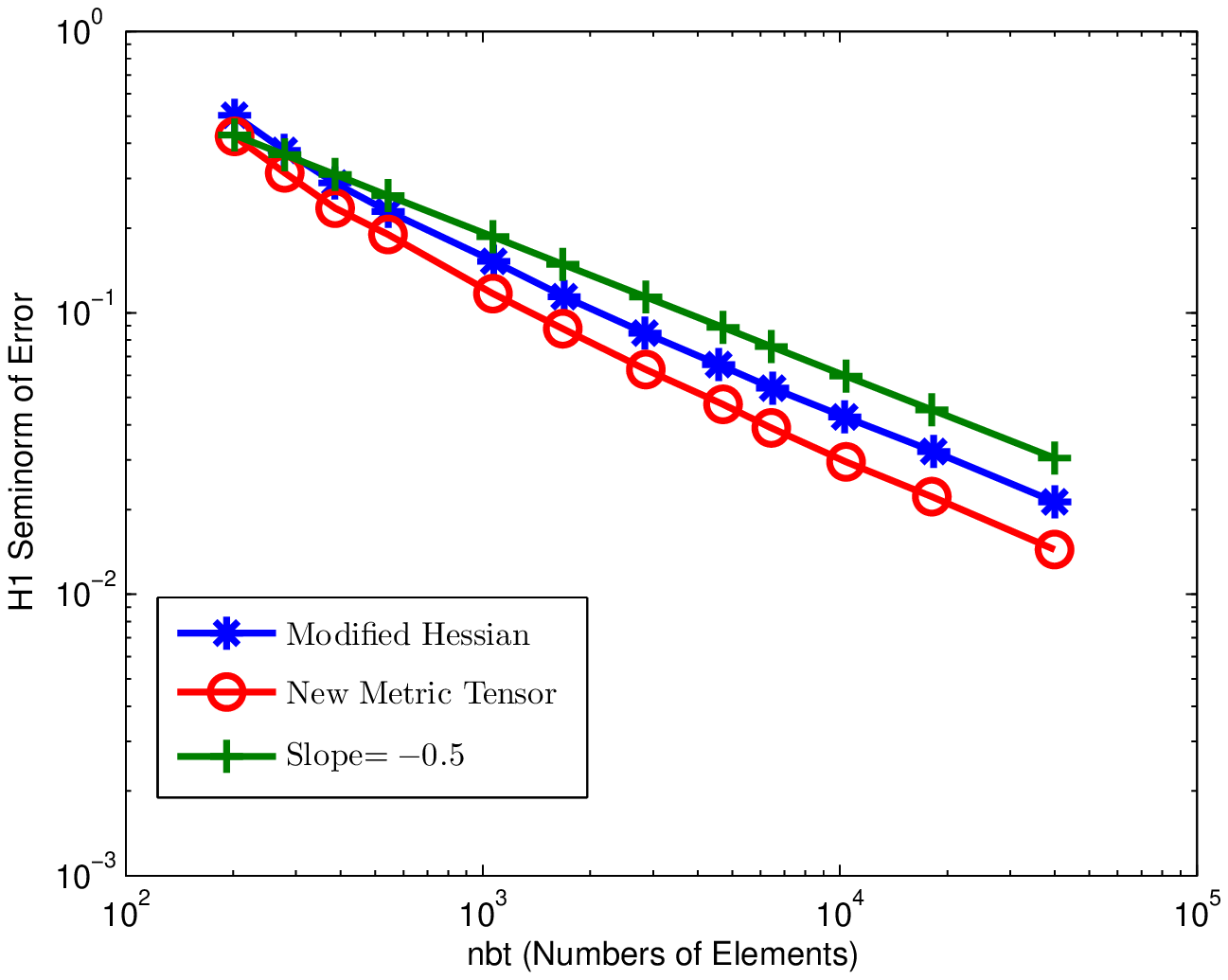}}}
    \put(-125,-15){(c)}
  \put(-20,-180){\resizebox{8cm}{!}{\epsffile{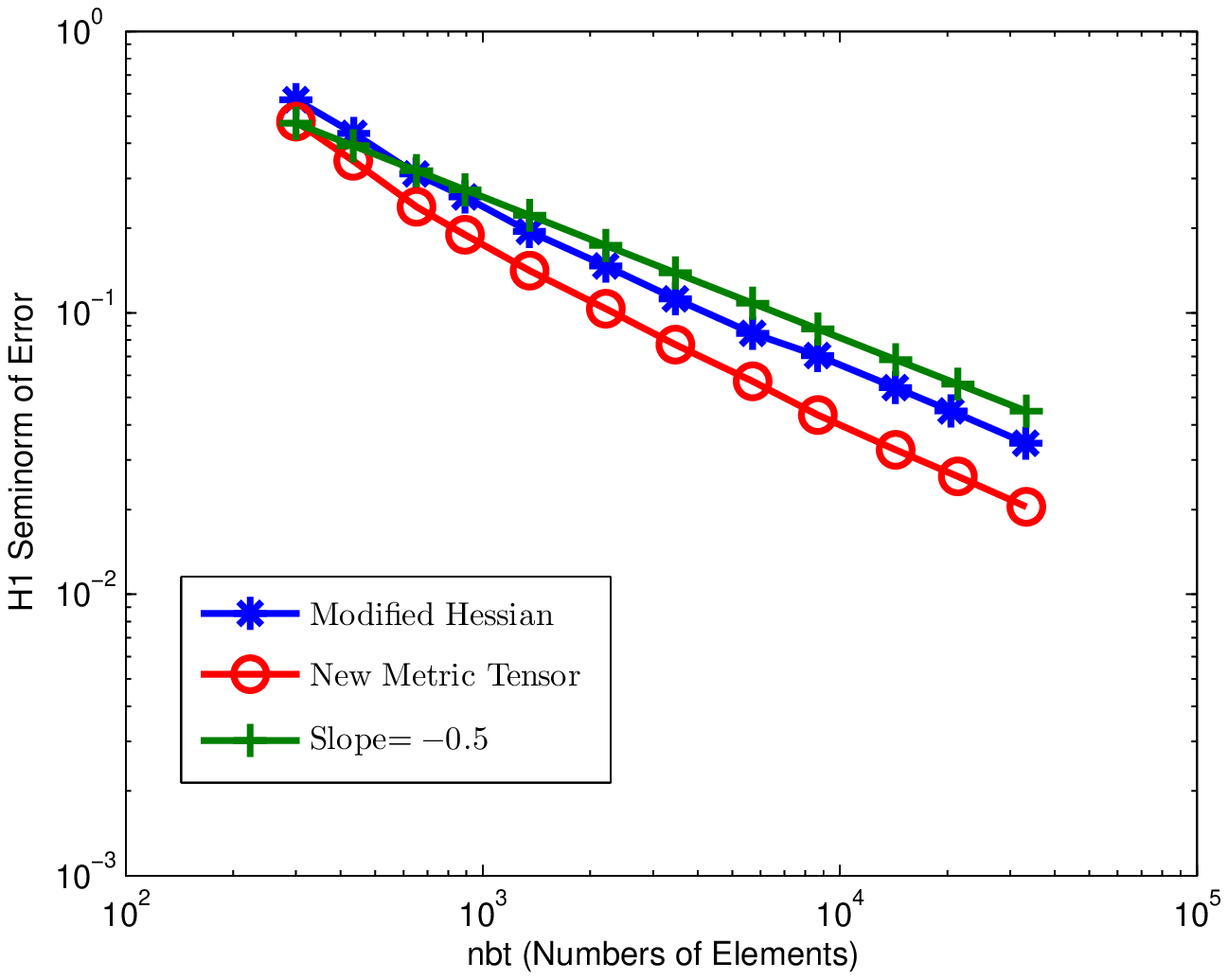}}}
    \put(100,-15){(d)}
    \caption{Example 4.3. The $H^1$ semi-norms of discretization error are plotted as functions
of $nbt$ under the conditions (a) $\beta=5$, (b) $\beta=10$, (c)
$\beta=20$, (d) $\beta=40$, respectively.}\label{Example3_H1error}
\end{figure}

\begin{figure}[ht]
 \includegraphics[width=8cm]{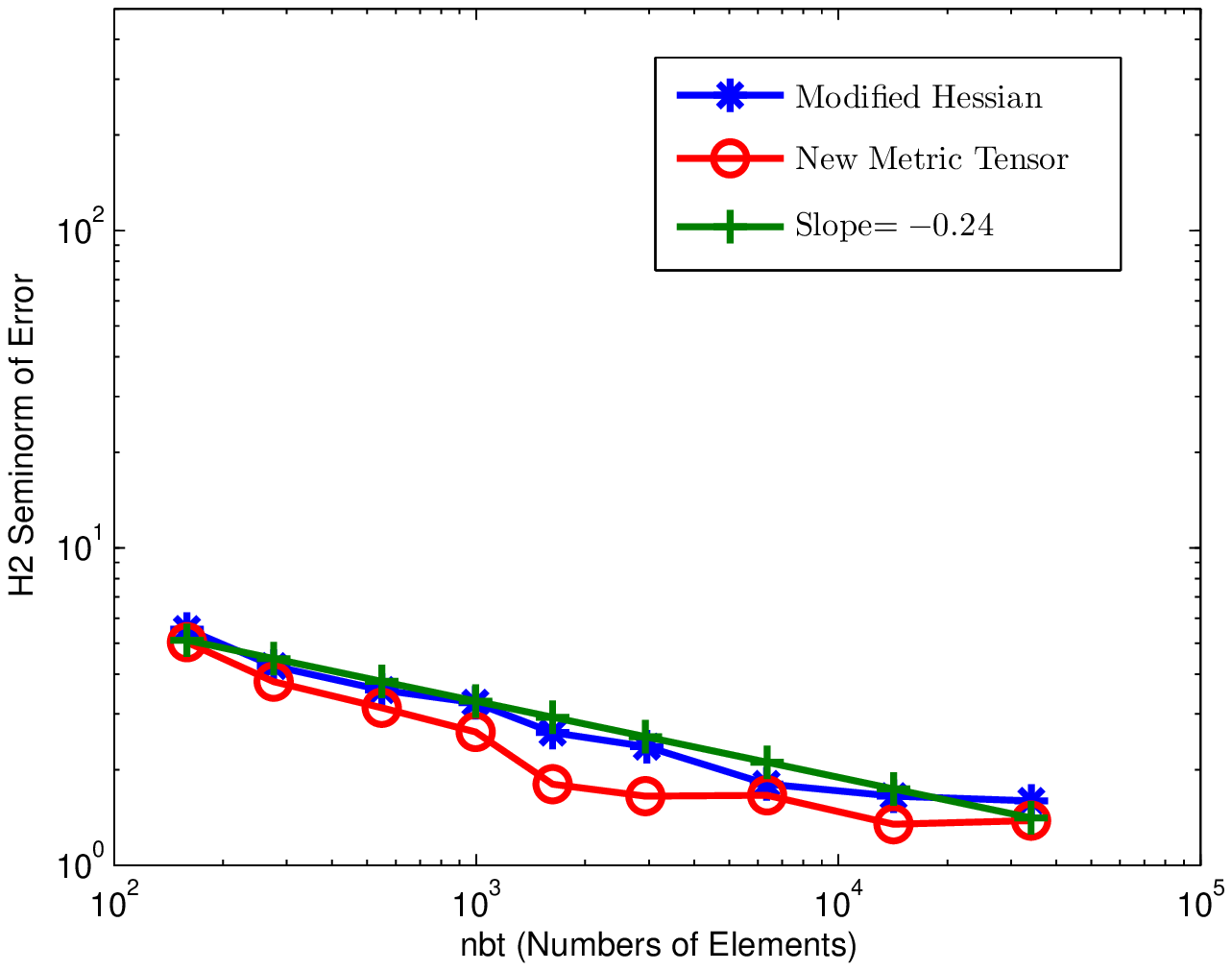}
  \put(-125,167){(a)}
  \put(-20,0){\resizebox{8cm}{!}{\epsffile{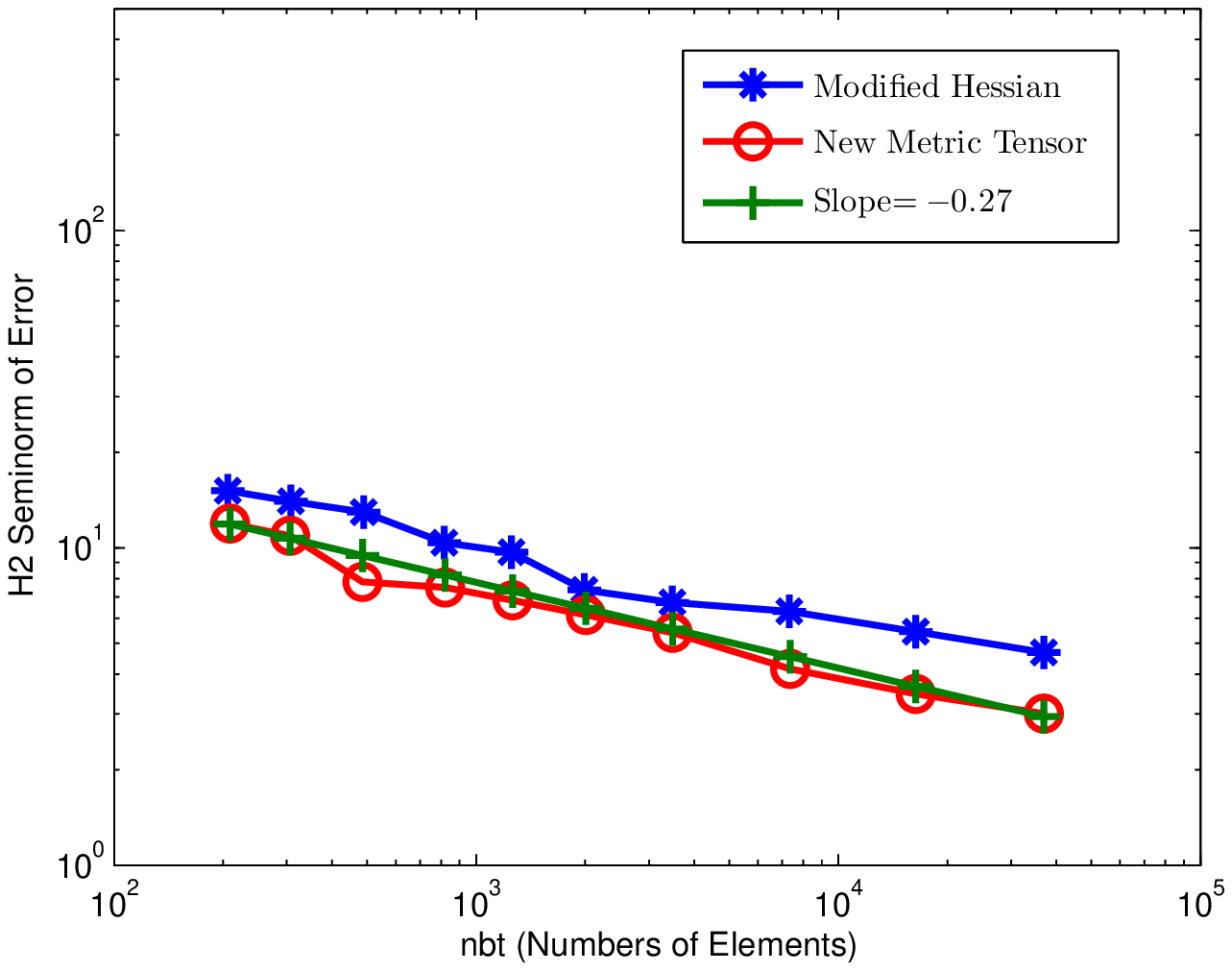}}}
  \put(100,167){(b)}
  \put(-228,-180){\resizebox{8cm}{!}{\epsffile{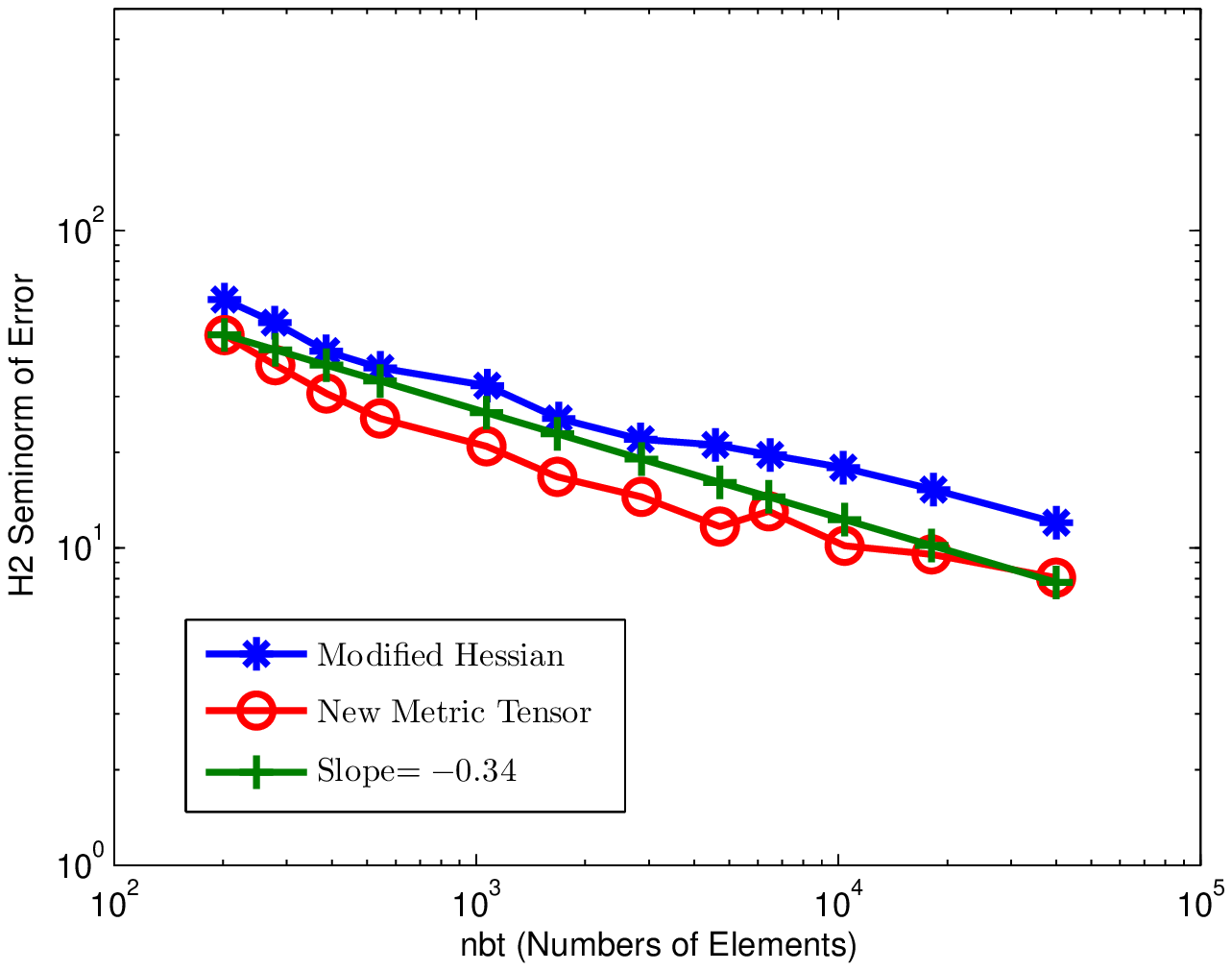}}}
  \put(-125,-15){(c)}
  \put(-20,-180){\resizebox{8cm}{!}{\epsffile{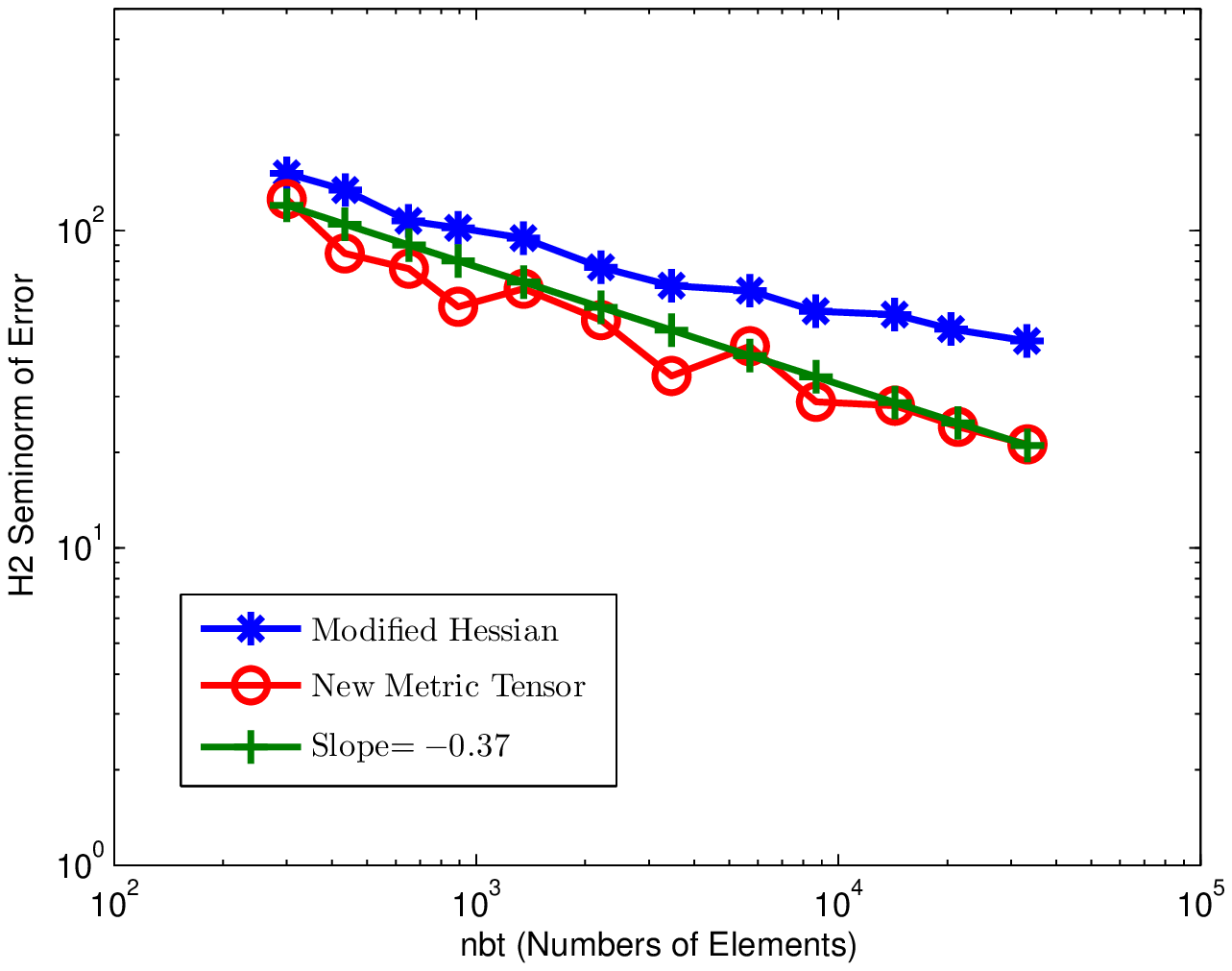}}}
   \put(100,-15){(d)}
   \caption{Example 4.3. The $H^2$ semi-norms of discretization error are plotted as functions
of $nbt$ under the conditions (a) $\beta=5$, (b) $\beta=10$, (c)
$\beta=20$, (d) $\beta=40$, respectively.}\label{Example3_H2error}
\end{figure}

Figure \ref{Example3_solution} contains two meshes obtained by the
two different metric tensors from solution perspective, where the
difference is obvious. This comparison, again, indicates our metric
tensor explores more comprehensive anisotropic information of the
solution $u$. Figure \ref{Example3_H1error} and Figure
\ref{Example3_H2error} contain $H^1$ and $H^2$ semi-norms of error,
respectively, under various parameters $\beta$. We know that the
larger the $\beta$ becomes, the more significant the anisotropy of
the solution exihibits. Then we can conclude that the more
anisotropic the solution is, the more obvious the difference is. In
fact, the difference lies on the term (\ref{factor}) which indicates
our metric tensor explores more comprehensive anisotropic
information of the solution when the term varies significantly at
different points or elements.

\section{Conclusions}
In the previous sections we have developed a new anisotropic mesh
adaptation strategy for finite element solution of elliptic
differential equations. Note that the new metric tensor
$\mathcal{M}_0$ for the $L^2$ norm is similar to (\ref{1.4_1})
(\cite{Huang}). However, the new metric tensor $\mathcal{M}_1$ is
essentially different with those metric tensors existed. The
difference lies on the term (\ref{1.5}) which indicates our metric
tensor explores more comprehensive anisotropic information of the
solution when the term varies significantly at different points or
elements. Numerical results also show that this approach is superior
to the existing ones to deal with poisson and steady
convection-dominated problems.


\begin{thebibliography}{23}
\bibitem{Ait-Ali-Yahia}
D. Ait-Ali-Yahia, W. Habashi, A. Tam, M.-G. Vallet, M. Fortin, A
directionally adaptive methodology using an edge-based error
estimate on quadrilateral grids, Int. J. Numer. Methods Fluids 23
(1996) 673-690.
\bibitem{ApelLube}
T. Apel, G. Lube, Anisotropic mesh
refinement in stabilized Galerkin methods, Numer. Math. 74(3) (1996)
261-282.
\bibitem{BankSmith}
R.E. Bank, R.K. Smith, Mesh smoothing using a posteriori error
estimates, SIAM J. Numer. Anal. 34 (1997) 979-997.
\bibitem{Becker}
R. Becker, An adaptive finite element method for the incompressible
Navier-stokes equations on time-dependent domains, Ph.D. thesis,
Ruprecht-Karls-Universit$\ddot{a}$t Heidelberg, 1995.
\bibitem{Borouchaki1}
H. Borouchaki, P.L. George, F. Hecht, P. Laug and E. Saltel,
Delaunay mesh generation governed by metric specifications Part I.
Algorithms, finite elem. anal. des. 25 (1997) 61-83.
\bibitem{Borouchaki2}
H. Borouchaki, P.L. George, B. Mohammadi, Delaunay mesh generation
governed by metric specifications Part II. Applications, finite
elem. anal. des. 25 (1997) 85-109.
\bibitem{Brackbill}
J.U. Brackbill, J.S. Saltzman, Adaptive zoning for singular problems
in two dimensions, J. Comput. Phys. 46 (1982) 342-368.
\bibitem{BuDa}
G. Buscaglia, E. Dari, Anisotropic Mesh Optimization and its
Application in Adaptivity, Int. J. Numer. Meth. Eng. 40(22) (1997)
4119-4136.
\bibitem{cao}
W. Cao, On the error of linear interpolation and the orientation,
aspect ratio, and internal angles of a triangle, SIAM J. Numer.
Anal. 43(1) (2005) 19-40.
\bibitem{CasHecMohPir}
M.J. Castro-D$\acute{\imath}$az, F. Hecht, B. Mohammadi, O.
Pironneau, Anisotropic unstructured mesh adaption for flow
simulations, Internat. J. Numer. Methods Fluids 25(4) (1997)
475-491.
\bibitem{ChenSunXu}
L. Chen, P. Sun, J. Xu, Optimal anisotropic meshes for minimizing
interpolation errors in $L\sp p$-norm, Math. Comp. 76(257) (2007)
179-204.
\bibitem{Azevedo1989}
E.F. D'Azevedo, R.B. Simpson, On optimal interpolation triangle
incidences, SIAM J. Sci. Statist. Comput. 10 (1989) 1063-1075.
\bibitem{Azevedo1991}
E.F. D'Azevedo, Optimal triangular mesh generation by coordinate
transformation, SIAM J. Sci. Stat. Comput. 12 (1991) 755-786.
\bibitem{AzevedoSimpson}
E.F. D'Azevedo, R.B. Simpson, On optimal triangular meshes for
minimizing the gradient error, Numer. Math. 59 (1991) 321-348.
\bibitem{Dvinsky}
A.S. Dvinsky, Adaptive grid generation from harmonic maps on
Riemannian manifolds, J. Comput. Phys. 95 (1991) 450-476.
\bibitem{ForPer}
L. Formaggia, S. Perotto,  Anisotropic error estimates for elliptic
problems, Numer. Math. 94(1) (2003) 67-92.
\bibitem{FreyGeorge}
P. Frey, P.L. George, Mesh Generation: Application to Finite
Elements, Hermes Science, Oxford and Paris, 2000.
\bibitem{Garimella}
R.V. Garimella, M.S. Shephard, Boundary layer meshing for viscous
flows in complex domain. in: Proceedings of the 7th International
Meshing Roundtable, Sandia National Laboratories, Albuquerque, NM,
1998, 107-118.
\bibitem{George}
P.L. George, F. Hecht. Nonisotropic grids, in: J.F. Thompson, B.K.
Soni, N.P. Weatherill, (Eds.), Handbook of Grid Generation, CRC
Press, Boca Raton, 1999 20.1-20.29.
\bibitem{HabDomBourAitForVal}
W.G. Habashi, J. Dompierre, Y. Bourgault, D. Ait-Ali-Yahia, M.
Fortin, M.-G. Vallet, Anisotropic mesh adaptation: towards
user-indepedent, mesh-independent and solver-independent CFD. Part
I: general principles, Int. J. Numer. Meth. Fluids 32 (2000)
725-744.
\bibitem{HabForDomValBou}
W.G. Habashi, M. Fortin, J. Dompierre, M.-G. Vallet, Y. Bourgault,
Anisotropic mesh adaptation: a step towards a mesh-independent and
user-independent CFD, Barriers and challenges in computational fluid
dynamics (Hampton, VA, 1996), 99-117, Kluwer Acad. Publ., Dordrecht,
1998.
\bibitem{Hecht}
F. Hecht, Bidimensional anisotropic mesh generator, Technical
Report, INRIA, Rocquencourt, 1997.
\bibitem{HuangSiam}
W. Huang. Measuring mesh qualities and application to variational
mesh adaptation. SIAM J. Sci. Comput. 26(5) (2005) 1643-1666.
\bibitem{Huang}
W. Huang,  Metric tensors for anisotropic mesh generation,  J.
Comput. Phys. 204(2) (2005) 633-665.
\bibitem{Jacquotte}
O.P. Jacquotte, A mechanical model for a new grid generation method
in computational fluid dynamics, Comput. Meth. Appl. Mech. Eng. 66
(1988) 323-338.
\bibitem{Knupp}
P. Knupp, L. Margolin, M. Shashkov, Reference jacobian
optimization-based rezone strategies for arbitrary lagrangian
eulerian methods, J. Comput. Phys. 176 (2002) 93-128.
\bibitem{Kornhuber}
R. Kornhuber, R. Roitzsch, On adaptive grid refinement in the
presence of internal or boundary layers, IMPACT Comput. Sci. Eng. 2
(1990) 40-72.
\bibitem{Lang}
J. Lang, An adaptive finite element method for convection-diffusion
problems by interpolation techniques, Technical Report TR 91-4,
Konrad-Zuse-Zentrum Berlin, 1991.
\bibitem{LiTangZhang}
R. Li, T. Tang, and P. Zhang, Moving mesh methods in multiple
dimensions based on harmonic maps, J. Comput. Phys. 170(2) (2001)
562-588.
\bibitem{Nadler}
E.J. Nadler, Piecewise linear approximation on triangulations of a
planar region, Ph.D. Thesis, Division of Applied Mathematics, Brown
University, Providence, RI, 1985.
\bibitem{PerVahMorZ}
J. Peraire, M. Vahdati, K. Morgan, O.C. Zienkiewicz, Adaptive
remeshing for compressible flow computation, J. Comp. Phys. 72(2)
(1987) 449-466.
\bibitem{Rachowicz}
W. Rachowicz, An anisotropic h-adaptive finite element method for
compressible Navier¨CStokes equations, Comput. Meth. Appl. Mech.
Eng. 146 (1997) 231-252.
\bibitem{Remacle}
J. Remacle, X. Li, M.S. Shephard, and J.E. Flaherty, Anisotropic
adaptive simulation of transient flows using discontinuous Galerkin
methods, Int. J. Numer. Meth. Eng., 62(7) (2005) 899-923.
\bibitem{Rippa}
S. Rippa, Long and thin triangles can be good for linear
interpolation, SIAM J. Numer. Anal. 29 (1992) 257-270.
\bibitem{Yamakawa}
S. Yamakawa and K. Shimada, High quality anisotropic tetrahedral
mesh generation via ellipsoidal bubble packing. in: Proceedings of
the 9th International Meshing Roundtable, Sandia National
Laboratories, Albuquerque, NM, 2000. Sandia Report 2000-2207,
263-273.
\bibitem{YinXie}
X. Yin, H. Xie, A-posteriori error estimators suitable for moving
mesh methods under anisotropic meshes, to appear.
\bibitem{Zhangxd}
X. Zhang, Accuracy concern for Hessian metric, Internal Note, CERCA.
\bibitem{Zhangzm}
Z. Zhang, A. Naga, A new finite element gradient recovery method:
superconvergence property, SIAM J. Sci. Comput. 26(4) (2005)
1192-1213.
\bibitem{ZWu}
O.C. Zienkiewicz, J. Wu, Automatic directional refinement in
adaptive analysis of compressible flows,  Int. J. Numer. Meth. Eng.
37 (1994) 2189-2210.
\bibitem{ZZ1987}
O.C. Zienkiewicz, J. Zhu, A simple error estimator and adaptive
procedure for practical engineering analysis, Int. J. Numer. Meth.
Eng. 24 (1987) 337-357.
\end{thebibliography}
\end{document}